\documentclass[a4paper, 10pt, final]{amsart}

\title[Algebraically bounded derivations]{Generic derivations on algebraically bounded structures}

\newcommand*{\mykeywords}{Derivation, algebraically bounded, model completion}

\usepackage{lmodern}

\author[A. Fornasiero]{Antongiulio Fornasiero}
\address{Universit\`a di Firenze}  
\email{antongiulio.fornasiero@gmail.com} 
\urladdr{https://sites.google.com/site/antongiuliofornasiero/}

\author[G. Terzo]{Giuseppina Terzo}
\address{Universit\`a degli Studi di Napoli "Federico II"}
\email{giuseppina.terzo@unina.it} 

\usepackage{ifpdf}
\ifpdf
\usepackage[pdftex,
pdfborder={0 0 0}, plainpages=false,%
pdfpagelabels, pdfdisplaydoctitle=true,%
pdfstartview={XYZ null null null},%
colorlinks=false,
pdfusetitle
]{hyperref}
\hypersetup{%
  pdfauthor={A. Fornasiero, G. Terzo},
  pdfkeywords={\mykeywords},
  pdfsubject={Derivations on an algebraically bounded structure,
model completion of the corresponding theory}
}
\else
\usepackage[plainpages=false,linktocpage=true,pdfpagelabels, colorlinks=false
]{hyperref}
\fi

\usepackage{amsmath, amsthm, amssymb}
\usepackage{mathtools}
\usepackage{xspace}

\usepackage{tikz-cd}

\usepackage{microtype}

\usepackage{comment}

\usepackage{enumitem}
\setlist{
         left = 0pt}
\setlist[enumerate,1]{label = (\arabic*), 
         ref =   (\arabic*),
         font = \upshape}

\newlist{customenum}{enumerate}{1}
\makeatletter
\def\ctext#1{\expandafter\@ctext\csname c@#1\endcsname}
\def\@ctext#1{\ifcase#1\or Deep\or Wide\or PP\or
PP-wrong\or Fifth\or Sixth\fi}
\makeatother
\AddEnumerateCounter{\ctext}{\@ctext}{Deep}
\setlist[customenum]{label=(\ctext*), 
  ref=(\texttt{\ctext*}),
  font= \texttt
}

\makeatother

\usepackage[T1]{fontenc}

\usepackage[english]{babel}
\usepackage[babel]{csquotes}
\usepackage[backend=biber, style=alphabetic, sorting=nyt,
 backref = true, backrefstyle = three,
 isbn = false, 
 ]{biblatex}
\AtEveryBibitem{%
  \clearlist{language}
  \clearfield{urlyear}
  \clearfield{urlmonth}
} 

\DeclareMathSymbol{\mlq}{\mathord}{operators}{``}
\DeclareMathSymbol{\mrq}{\mathord}{operators}{`'}

\DeclareMathOperator{\Reg}{Reg}

\newcommand{\thefree}{\mathcal F}
\newcommand{\thebounded}{\mathcal B}
\newcommand{\thepivot}{\mathcal P}
\newcommand{\tildefree}{\mathcal V}
\newcommand{\posture}{configuration\xspace}
\newcommand{\postures}{configurations\xspace}
\DeclareMathOperator{\DCFO}{DCF_0}
\DeclareMathOperator{\DCFOm}{DCF_{0,m}}
\DeclareMathOperator{\DCFOmnc}{DCF_{0,m,nc}}

\DeclareMathOperator{\Diag}{Diag}

\newcommand{\K}{\mathbb K}

\newcommand{\Td}{T^{\delta}}
\newcommand{\Tdg}{T^{\delta}_g}
\newcommand{\Tdgp}{{T^{\delta}_{\mathrm{deep}}}}
\newcommand{\Tdgpp}{{T^{\delta}_{\mathrm{wide}}}}
\newcommand{\TdgPP}{{T^{\delta}_{\mathrm{PP}}}}
\newcommand{\Tdb}{T^{\bar \delta}}
\newcommand{\Tdbg}{T^{\bar \delta}_g}
\newcommand{\Tdn}{T^{\bar \delta,nc}}
\newcommand{\Tdng}{T^{\bar \delta,nc}_g}
\newcommand{\Tdngp}{T^{\bar \delta,nc}_{deep}}
\newcommand{\Tdngpp}{T^{\bar \delta,nc}_{wide}}

\newcommand{\Tdsg}{T^{\bar \delta,?}_g}

\newcommand{\condition}{\mathfrak S}

\newcommand{\Ldb}{L^{\bar \delta}}
\newcommand{\deltab}{\bar\delta}
\newcommand{\epsb}{\bar\eps}

\newcommand{\eps}{\varepsilon}
\DeclareMathOperator{\rk}{rk}
\DeclareMathOperator{\Jet}{Jet}

\newcommand{\Jetdn}{\Jet^n_{\delta}}
\newcommand{\Jetdinf}{\Jet^{\infty}_{\delta}}

\newcommand{\Ldelta}{L^\delta}
\newcommand{\Ldeltabar}{L^{\bar \delta}}
\newcommand{\monster}{\mathbb M}

\newcommand*{\intro}[1]{\textbf{#1}}
\newcommand*{\Pa}[1]{\bigl( #1 \bigr)}
\newcommand*{\set}[1]{\{#1\}}
\newcommand*{\abs}[1]{\lvert#1\rvert}

\newcommand*{\card}[1]{\lvert#1\rvert}
\newcommand{\N}{\mathbb{N}}

\newcommand{\Q}{\mathbb{Q}}

\newcommand{\rest}{\upharpoonright}

\DeclareMathOperator{\dcl}{dcl}
\DeclareMathOperator{\acl}{acl}
\newcommand*{\tuple}[1]{\langle #1 \rangle}

\newcommand{\av}{\bar a}
\newcommand{\bv}{\bar b}
\newcommand{\cv}{\bar c}
\newcommand{\dv}{\bar d}

\newcommand{\pv}{\bar p}
\newcommand{\x}{\bar x}
\newcommand{\y}{\bar y}

\newcommand{\deltabar}{\bar \delta}

\def\Ind#1#2{#1\setbox0=\hbox{$#1x$}\kern\wd0\hbox to
  0pt{\hss$#1\mid$\hss}\lower.9\ht0\hbox to 0pt{\hss$#1\smile$\hss}\kern\wd0}

\newcommand{\functiondef}[5]{%
\begin{align*}
  {#1} : {#2} & \longrightarrow  {#3} \\
  {#4} &\longmapsto\! {#5}
\end{align*}
}

\newcommand{\et}{\ \wedge\ }
\newcommand{\vel}{\ \vee\ }

\def\hyph{\nobreakdash-\hspace{0pt}\relax}

\newcommand{\Wlog}{W.l.o.g\mbox{.}\xspace}
\newcommand{\wloG}{w.l.o.g\mbox{.}\xspace}

\newcommand{\ie}{i.e\mbox{.}\xspace}

\newcommand{\wrt}{w.r.t\mbox{.}\xspace}

\newtheorem{lemma}{Lemma}[section]
\newtheorem{thm}[lemma]{Theorem}
\newtheorem{corollary}[lemma]{Corollary}
\newtheorem{conjecture}[lemma]{Conjecture}
\newtheorem{proposition}[lemma]{Proposition}
\newtheorem{open problem}[lemma]{Open problem}

\newtheorem*{fact*}{Fact}
\newtheorem{fact}[lemma]{Fact}
\newtheorem{def-lemma}[lemma]{Definition-Lemma}

\theoremstyle{remark}

\newtheorem{claim}{Claim}
\newtheorem*{claim*}{Claim}

\theoremstyle{definition}
\newtheorem{definition}[lemma]{Definition}

\newtheorem{remark}[lemma]{Remark}
\newtheorem{final remark}[lemma]{Final remark}

\newtheorem{question}[lemma]{Question}

\newenvironment{sentence}[1][]{%
  \begin{list}{}{%
    \setlength\topsep{0.5ex}%
    \setlength\leftmargin{\parindent}%
  }%
  \item[#1]
 }
 {\end{list}}

\begin{document}

\begin{abstract}
Let $\K$ be an algebraically bounded structure and $T$ be its theory.
If $T$ is model complete, then the theory of $\K$ endowed with a derivation, denoted by  $\Td$, has a
model completion. Additionally, we prove that if the theory $T$ is stable/NIP then the model completion  of $\Td$ is also stable/NIP. Similar results hold for the theory with several  derivations, either commuting or non-commuting.


\end{abstract}



\maketitle


\makeatletter
\renewcommand\@makefnmark%
   {\normalfont(\@textsuperscript{\normalfont\@thefnmark})}
\renewcommand\@makefntext[1]%
   {\noindent\makebox[1.8em][r]{\@makefnmark\ }#1}
\makeatother

\tableofcontents

\section{Introduction}

Let $\K$ be a structure expanding a field of characteristic 0.
Recall that $\K$ is \textbf{algebraically bounded} if the model-theoretic algebraic
closure and the field-theoretic algebraic closure coincide in every structure elementarily equivalent to~$\K$.
Algebraically closed, real closed, p-adically closed, pseudo-finite fields, and
algebraically closed valued fields
are examples of algebraically bounded structures; for more details, examples, and main properties
see~\cite{Dries:89} and \S\ref{sec:AB}. Van den Dries in his paper introduced a notion of dimension for any definable set with parameters, which is relevant in our context.  

Let $L$ be the language of $\K$ and $T$ be its theory. 
In order to study derivations on~$\K$, we denote by $\delta$ a new unary function symbol, 
and by $\Td$ the $\Ldelta$-theory expanding $T$ by saying that $\delta$ 
is a derivation.  
Let $\K$ be algebraically bounded. We remark that being algebraically bounded is not a first order notion (since an ultraproduct of algebraically bounded structures is not necessarily algebraically bounded). We define an $\Ldelta$-theory $\Tdg$
extending $\Td$, with three equivalent axiomatizations
(see \S\S\ref{sec:MC}, \ref{sec:PP}); one of them is given by $\Td$, plus the 
following axiom scheme:\\
\begin{sentence}
For every $X \subseteq \K^{n} \times \K^{n}$ which is $L$-definable with parameters, if the
dimension of the projection of $X$ onto the first $n$ coordinates, which we denote by $\Pi_{n}(X),$  is~$n$,  then there
exists $\av \in \K^{n}$ such that $\tuple{\av, \delta \av} \in X$.\\
\end{sentence}

One of the main result of the paper is the following:

\begin{thm}
 If $T$ is model complete, then $\Tdg$ is the model completion of~$\Td$.
 \end{thm}


We also endow $\K$ with several derivations $\delta_{1}, \dotsc , \delta_{m}$
and we  consider both the case when they commute and
when we don't impose any commutativity.
We obtain two theories:
\begin{description}
\item[$\Tdb$] the expansion of $T$ saying that the $\delta_{i}$ are derivations which
commute with each other;
\item[$\Tdn$] the expansion of $T$ saying that the
$\delta_{i}$ are derivations without any further conditions.
\end{description}
Both theories have a model completion (if $T$ is model complete) 
(see \S\S\ref{sec:commute}, \ref{sec:nc}).
For convenience, we use $\Tdsg$ to denote either of the model completions, both
for commuting derivations and the non-commuting case.
Many of the model-theoretic properties of $T$ are inherited by~$\Tdsg$:
\begin{thm}[\S\ref{sec:stable}] 
Assume that $T$ is stable/NIP.
Then  $\Tdsg$ is stable/NIP. 
\end{thm}

In a work in preparation, we will prove that if T is simple, then $\Tdsg$ is
simple (see \cites{MS, Mohamed} for particular cases); we will also characterize when  $\Tdsg$ is $\omega$-stable.
Moreover, we will prove that $\Tdsg$ is uniformly finite. 
Finally, if $\K$ extends either a Henselian valued field with a definable
valuation or a real closed field, then, under some additional assumptions,  
$T$~is the open core of $\Tdsg$.

On the other hand, in \cite{FT:dexp} we show that exponential fields do not
admit generic derivations (notice that exponential fields are not algebraically bounded).

\medskip

\subsection{A brief model theoretic history}\label{storico}

From a model theoretic point of view, differential fields have been studied
at least since Robinson \cite{Robinson} proved that the theory of fields of characteristic 0
with one derivation has a model completion, the theory $\DCFO$ of differentially
closed fields of characteristic~$0$. 

Blum gave a simpler sets of axioms for $\DCFO$, saying that $\K$ is
a  field of characteristic $0$, and, whenever $p$ and $q$ are
differential polynomials in one variable, with $q$ not zero and of order
strictly less than the order of~$p$, then there exists $a$ in $\K$
such that $p(a) = 0$ and $q(a) \neq 0$ (see \cite{Blum:77, Sacks} for more details).
Pierce and Pillay \cite{PP} gave yet another axiomatization for $\DCFO$, which has been
influential in the axiomatizations of other structures  (see \S\ref{sec:PP}).

The theory $\DCFO$ (and its models) has been studied intensively, both for its
own sake, for applications, and
as an important example of many ``abstract'' model theoretic properties: it is $\omega$-stable of rank $\omega$, it eliminates
imaginaries, it is uniformly finite, etc.
For some surveys see \cites{Blum:77, HV,  Chatzidakis, MMP, Moosa:22}.

Models of $\DCFO$, as fields, are algebraically closed fields of
characteristic $0$;
their study has been extended in several directions.
An important extension, which however goes beyond the scope of this article, is
Wood's work \cite{Wood}  on fields of finite characteristic. 
\medskip

From now on all
fields are of characteristic~$0$. More close to the goal of this article is the passage from one derivation to
several commuting ones: McGrail  \cite{mcgrail} axiomatized $\DCFOm$ (the model completion of the theory of fields of characteristic 0 with $m$ commuting derivations). 
While the axiomatization is complicate (see \S\ref{sec:commute} for an easier
axiomatization, and \cites{Pierce14, LS} for alternative ones), from a model
theoretic point of view $\DCFOm$ is quite similar to $\DCFO$: its models
are algebraically closed (as fields), it is  $\omega$-stable of rank
$\omega^{m}$, it eliminates imaginaries, it is uniformly finite, etc.

Moosa and Scanlon followed a different path in \cite{MS}, where they studied a
general framework of fields with non-commuting opeators; 
for this introduction, the relevant
application is that they proved that  the theory of  $m$
non-commuting derivations has a model completion (see \cite{MS} and
\S\ref{sec:nc}), which we denote by $\DCFOmnc$.
Here the model theory is more complicate: $\DCFOmnc$ is stable,  
but not $\omega$-stable; however, it still eliminates
imaginaries and it is uniformly finite.\\
Surprisingly, we can give 3 axiomatizations for $\DCFOmnc$ which are much simpler
than the known axiomatizations for $\DCFOm$ (including the one given in
this article), see \S\S\ref{sec:nc}, \ref{sec:PP}.
We guess that the reason why this has not been observed before is that
people were deceived by the rich algebraic structure of $\DCFOm$.\\
Indeed, from an algebraic point of view, $\DCFOm$ 
has been studied extensively (see
\cite{Kolchin} for a starting point) and is much simpler than $\DCFOmnc.$ 
The underlying combinatorial fact is that the free commutative monoid on $m$
generators $\Theta$, with the partial ordering given by
$\alpha \preceq \beta \alpha $ for every $\alpha, \beta \in \Theta$, is a well-partial-order (by Dickson's Lemma);
this fact is a fundamental ingredient in Ritt-Raudenbush Theorem, asserting that
there is no infinite ascending chain of radical differential ideals in the ring
of differential polynomials with $m$ commuting derivations with coefficients in
some differential field; 
moreover, every radical differential ideal is a finite intersection
of prime differential ideals.  Since in models of $\DCFOm$ there is a natural
bijection between prime differential ideals and complete types, this in turns
implies that $\DCFOm$ is $\omega$-stable as we mentioned before.

Very different is the situation for the free monoid on $m$ generators $\Gamma$,
with the same partial ordering.
$\Gamma$ is well-founded, but (when $m$ is at least 2) not a well-partial-order.
Given an infinite anti-chain in~$\Gamma$, it is easy to build an infinite ascending
chain of radical differential ideals (in the corresponding ring
of non-commuting differential polynomials),
and therefore Ritt-Raudenbush does not hold in this situation.

Some limited form of non-commutativity was considered already in 
\cites{Yaffe:01, Singer:07, Pierce14}, 
where the derivations live in a finite-dimensional Lie algebra.

\medskip

People have extended $\DCFO$ in another direction by considering fields which are
not algebraically closed: Singer, and later others
\cites{Singer:78, Point, bkp, BMR, Riviere:09}  studied real closed fields with
one generic derivation, and \cite{Riviere:06b} extended to $m$ commuting
derivations (see also \cite{FK} for a different approach); 
\cites{GP:10,GP:12,GR,KP} studied more general topological fields
with one generic derivation. In \cite{Riviere:06} the author studied fields with $m$ independent
orderings and one generic derivation and in 
\cite{FK} they studied o-minimal structures with several commuting generic 
``compatible'' derivations.
In her Ph.D.\  thesis, Borrata \cite{Borrata:phd} studied ordered valued fields and
``tame'' pairs of real closed fields endowed with one generic derivation.

The results in \cite{KP,GP:10,GP:12,KP,Riviere:06, Riviere:06b} extend the one in \cite{Singer:78} and are mostly subsumed in  this article
(because the structures they study are mostly algebraically bounded).

\medskip

Tressl in \cite{Tressl:05} studied generic derivations on fields that are
``large'' in the sense of Pop, and Mohamed in \cite{Mohamed} extended his work to
operators in the sense of~\cite{MS}. 
They assume that their fields are model-complete in the language of rings with
additional constants, and therefore they are algebraically bounded
(see \cite[Thm.5.4]{JK}).%
\footnote{There is a slight misstatement in their theorem, in
that $\K$ must be in the language  rings with constants, and not only
a ``pure'' field as defined in their paper; besides, their proof allows adding
constants to the language in characteristic~0).}
Thus, our results extend their result on the existence of a model companion for
large fields with finitely many derivations (either commuting as in
\cite{Tressl:05} or non-commuting as in~\cite{Mohamed}). 
Moreover, in this paper we consider fields which are not
pure fields, such as algebraically closed valued fields (see \S\ref{sec:AB} for
more examples).


It turns out that, while in practice many of the fields studied in model theory
are both large and algebraically bounded (and therefore their generic
derivations can be studied by using either our framework or the one of
Tressl et al.), there exist large fields which are not algebraically bounded
(the field $\mathbb C((X,Y))$ is large but not algebraically bounded, 
see \cite[Example~8]{Fehm}),
and there exist algebraically bounded fields which are not large 
(see~\cite{JY:geometric}).
Tressl and Le\`on Sanchez  \cites{ST:24, ST:23}
later introduce the notion of ``differentially large fields''.

\smallskip

Often the fields considered have a topology (e.g. they are ordered fields or valued fields):
however, the theories described above do not impose any continuity on the derivation (and the corresponding  ``generic'' derivations are not continuous at any point).
In \cite{Scanlon1, Scanlon2} and \cite{ADH} the authors consider the case of a valued field endowed with a  ``monotone'' derivation (i.e. a derivation $\delta$ such that $v(\delta x) \geq v(x);$ in particular, $\delta$ is continuous) and prove a corresponding Ax-Kochen-Ersov principle.


\section{Algebraically boundedness and dimension}\label{sec:AB}

We fix an L-structure $\K$ expanding a field of characteristic $0$.

We recall the following definition in  \cite{Dries:89}, as refined in \cite{JY:geometric}:

\begin{definition}
Let $F$ be a subring of $\K$.
We say that $\K$ is algebraically bounded over F if, for any formula $\phi(\x, y)$, there exist finitely many polynomials $p_1, \ldots, p_m \in F[\x, y]$ such that for any $\av,$ if $\phi(\av, \K)$ is finite, then $\phi(\av, \K)$ is contained in the zero set of $p_i(\av, y)$ for some $i$ such that $p_i(\av, y)$ doesn't vanish. $\K$ is \textbf{algebraically bounded} if it is algebraically bounded over $\K$.
\end{definition}

Since we assumed that $\K$ has characteristic~$0$, 
in the above definition 
we can replace ``$p_i(\av, y)$ doesn't vanish'' with the following:\\
``$p_i(\av, b) = 0$ and
$\frac{\partial p_{i}}{\partial y}(\av, b) \neq 0$''.

\begin{fact}[\cite{JY:geometric}, see also \cite{Fornasiero:matroid}]
\begin{enumerate}
\item $\K$ is algebraically bounded iff it is algebraically bounded over
$\dcl(\emptyset)$.
\item $\K$ is algebraically bounded
over~$F$ iff the model theoretic algebraic closure coincide with the field theoretic
algebraic closure over $F$ in every  elementary extension of $\K$
(it suffices to check it in the monster model). 
\end{enumerate}
\end{fact}

\begin{remark}
Junker and Koenigsmann in \cite{JK} defined $\K$ to be ``very slim'' if in the monster model
the field-theoretic algebraic closure over the prime field coincide with the
model-theoretic algebraic closure: thus, $\K$ is very slim iff $\K$ is
algebraically bounded over~$\Q$.
\end{remark}

\smallskip
Let $F \coloneqq \dcl(\emptyset)$ and we consider $\K$ algebraically bounded over~$F$.

\smallskip

When we  refer to the algebraic closure, unless specified otherwise, we will
mean the $T$-algebraic closure; similarly, $\acl$ will be the $T$-algebraic
closure, and by ``algebraically independent'' we will mean according to $T$ (or
equivalently algebraically independent over $F$ in the field-theoretic meaning).

From the assumptions it follows that $\K$ is \intro{geometric}: that is, in the
monster model $\monster \succ \K$, 
the algebraic closure has the exchange property, and therefore it is a
matroid; moreover, $T$ is Uniformly Finite, that is it eliminates the
quantifier $\exists^{\infty}$. 
In fact, a definable set $X \subseteq \K$ is infinite iff for all $a \in K$ there exist $x, y, x', y' \in X$ such that $x \not = x'$  and $a = (y-y')/(x'-x);$ 
(see  \cites{JY:geometric, Fornasiero:matroid}).

Moreover, $\K$ is endowed with a dimension function $\dim,$ associating to every
set $X$ definable with parameters some natural number, satisfying the axioms in
\cite{Dries:89}.
This function $\dim$ is invariant under automorphisms of the ambient
structure: equivalently, $\dim$ is ``code-definable'' in the sense of~\cite{bkp}.

We will also use the rank, denoted by $\rk,$ associated to the matroid $\acl$: $\rk(V/B)$ is
the cardinality of a basis of $V$ over $B$.
Thus, if $X \subseteq \monster^{n}$ is definable with parameters~$\bv$,
\[
\dim(X) = \max \Pa{\rk(\av/ \bv): \av \in X}.
\]

\subsection{Examples}

Some well known examples of fields which are  algebraically bounded structures as
pure fields are: algebraically closed fields, $p$-adics and more generally
Henselian fields (see \cite[Thm~5.5]{JK}), real closed fields,
pseudo-finite fields; curve-excluding fields in the sense of \cite{Johnson1} are
also algebraically bounded.
Other examples of algebraically bounded structures which are not necessarily
pure fields are:
\begin{itemize}
\item Algebraically closed valued fields; 
\item Henselian fields (of characteristic 0)
with arbitrary relations on the value group and the residue field (see
\cite{Dries:89});
\item All models of ``open theories of topological fields'', as defined in \cite{KP};
\item The expansion of an algebraically bounded structure by a
generic predicate (in the sense of \cite{CP:98}) is still algebraically bounded
(see \cite[Corollary~2.6]{CP:98});
\item The theory of fields with several independent orderings and valuations has
a model companion, whose models are algebraically bounded (see \cite{Dries:phd}, 
\cite[Corollary~3.12]{Johnson:22}).
\end{itemize}

Johnson and Ye in a recent paper
\cite{Johnson1} 
produced examples of an infinite algebraically bounded field with a
decidable first-order theory which is not large (in the Pop sense), and of a pure
field that is algebraically bounded but not very slim.

\subsection{Assumptions}


Our assumptions for the whole article are the following:

\begin{itemize}
\item $\K$ is a structure expanding a field of characteristic $0$.  
\item $L$ is the language
of $\K$ and $T$ is its $L$-theory.
\item $F \coloneqq \dcl(\emptyset) \subseteq \K$.
\item $\K$ is algebraically bounded over~$F$.  
\item $\dim$ is the dimension function on $\K$ (or on any model of $T$), $\acl$
is the $T$-algebraic closure, and $\rk$ the rank of the corresponding matroid.
\end{itemize}

\section{Generic derivation}\label{sec:MC}

We  fix a derivation $\eta: F \to F$ (if $F$ is contained in the algebraic closure of $\Q$ in~$\K$, that we denote by $\overline \Q,$ then $\eta$ must be equal to $0$).
We denote by $\Td$ the expansion of~$T$,
saying that $\delta$ is a derivation on $\K$ extending~$\eta$.

In the most important case, $F = \overline \Q$ and therefore $\eta = 0$, and $\Td$ is the
expansion of $T$ saying that $\delta$ is a derivation on~$\K$.

\subsection{ Model completion}
A. Robinson introduced the notion of model completion in relation with solvability of systems of equations. For convenience  we recall the definition:
\begin{definition} \label{modelcompletion}
Let $U$ and $U^{*}$ be theories in the same language~$L.$\\
$U^{*}$ is a model completion of $U$ if the following hold:
\begin{enumerate}
\item If $A \models U^{*},$ then $A \models U;$
\item If $A \models U,$ then there exists a $B \supset A$ such that $B \models U^{*};$
\item If $A \models U,$ $A \subset B,$ $A \subset C,$ where $B, C  \models U^{*},$ then B is elementary equivalent to  C over A.
\end{enumerate}
\end{definition}

We give the following general criteria for model completion. In our context we use \ref{en:Blum}.

\begin{proposition}\label{prop:MC}
Let $U$ and $U^{*}$ be theories in the same language~$L$ such that
$U \subseteq U^{*}$.
The following are equivalent:
\begin{enumerate}
\item\label{eq} $U^{*}$ is the model completion of~$U$ and $U^{*}$ eliminates quantifiers.

\item\label{blum1}
\begin{enumerate}
\item \label{blum1a}
For every $A \models U$, for every
$\sigma_{1}, \dotsc, \sigma_{n} \in U^{*}$, there exists $B \models U$ such that
$A \subseteq B$ and $B \models \sigma_{1} \wedge \dots  \wedge  \sigma_{n}$;
\item  \label{blum1b} For every $L$-structures $A, B, C$ such that $B \models U$,
$C \models U^{*}$, and $A$ is a common substructure, for every quantifier-free
$L(A)$-formula $\phi(\x)$, for every $\bv \in B^n$ such that
$B \models \phi(\bv)$, there exists $\cv \in C^n$ such that $C \models \phi(\cv)$.
\end{enumerate}

\item\label{en:Blum}
\begin{enumerate}
\item\label{en:Bluma} For every $A \models U$, for every
$\sigma_{1}, \dotsc, \sigma_{n} \in U^{*}$, there exists $B \models U$ such that
$A \subseteq B$ and $B \models \sigma_{1} \wedge \dots \wedge \sigma_{n}$;
\item \label{en:Blumb}For every $L$-structures $A, B, C$ such that $B \models U$,
$C \models U^{*}$, and $A$ is a common substructure, for every quantifier-free
$L(A)$-formula $\phi(x)$, and for every $b \in B$ such that $B \models \phi(b)$, there exists
$c \in C$ such that $C \models \phi(c)$.
\end{enumerate}
\item\label{diagramma}For all models A of $U_{\forall}$ we have:
\begin{enumerate}
\item  $\Diag(A) \cup U^{*}$ is
consistent,
\item  $\Diag(A) \cup U^{*}$ is
complete,
\end{enumerate}
where $\Diag(A)$ is the $L$-diagram of $A.$
\item\label{blum} (Blum criterion) \begin{enumerate}

\item Any model of $U_{\forall}$  can be extended to some model of $U^{*}$ . 
\item For any $A, A(b) \models U_{\forall}$ and for all $C^{*} \models U^{*},$ where $C^{*}$ is $|A|^+$-saturated, there exists an immersion of $A(b)$ in  $C^{*}$. 
 
\end{enumerate}

\item\label{universale} $U^{*}$ is the model completion of $U_{\forall}$.

\end{enumerate}
\end{proposition}

\begin{proof}
First we prove that \ref{eq} is equivalent to \ref{universale}: if $U^{*}$ is
the model completion of $U_{\forall},$ trivially $U^{*}$ is a model completion of $U$
and by \cite[Thm.~13.2]{Sacks}, we have that $U^{*}$ eliminates quantifiers.\\
For the converse, we have trivially that any models of $U^{*}$ is a model 
of~$U_{\forall}$. 
Moreover, if $A \models U_{\forall}$ then there exists $C \models U$ such that there
exists an immersion of $A$ in~$C$. 
But by \ref{universale} there exists $B \models U_{\forall}$ such that there exists an
immersion of $C$ in $B,$ and so an immersion of $A$ in $B$. 
It is trivial to verify (3) in Definition~\ref{modelcompletion}. 
\ref{eq} is equivalent to \ref{diagramma}: see~\cite{Sacks}. 
Also for the equivalence between \ref{blum} and \ref{universale} 
see~\cite{Sacks}.\\
It remains to prove the equivalence between \ref{eq} and
\ref{blum1}. 
$\ref{eq} \Rightarrow \ref{blum1} \Rightarrow \ref{en:Blum}$ is easy. 
For $\ref{en:Blum} \Rightarrow \ref{eq}$, in order to obtain that $U^{*}$ is the model
completion of~$U$ we prove that $\Diag(A) \cup U^{*}$ is consistent, but it is
enough to see that it is finitely consistent. By \ref{en:Bluma} we have the
finitely consistent. To prove that $U^{*}$ eliminates quantifiers it is
equivalent to prove that $\Diag(A) \cup U^{*}$ is complete, which follows easily
from~\ref{en:Blumb}.
\end{proof}

\subsection{The axioms}
We introduce the following notation:\\
Let $\delta: \K \to \K$ be some function, $n \in \N$, $a \in \K$
and $\av$ tuple of $\K^n.$ We denote by

\[
\begin{aligned}
\Jetdinf(a) & \coloneqq \tuple{\delta^{i} a: i \in \N}, &
\Jetdn(a) & \coloneqq \tuple{\delta^{i} a: i \leq n}, & \Jet(a) := \Jetdn(a) \mbox{ for some } n,\\
\Jetdinf(\av) & \coloneqq \tuple{\delta^{i} \av: i \in \N}, &
\Jetdn(\av) & \coloneqq \tuple{\delta^{i} \av: i \leq n}, & \Jet(\av) := \Jetdn(\av) \mbox{ for some } n.
\end{aligned}\\
\]

\begin{definition}
Let $X \subseteq \K^{n}$ be $L$-definable with parameters.
We say that $X$ is \intro{large} if $\dim(X) = n$.
\end{definition}

Two possible axiomatizations for the model completion $\Tdg$ are given by $\Td$ and either of the
following axiom schemas:
\begin{customenum}
\item\label{deep}
For every $Z \subseteq \K^{n+1}$ $L(\K)$-definable, if $\Pi_{n}(Z)$ is large, then there exists
$c \in \K$ such that $\Jetdn(c) \in Z$;
\item\label{wide}
For every $W \subseteq \K^{n} \times \K^{n}$ $L(\K)$-definable, if $\Pi_{n}(W)$ is large, then there exists
$\cv \in \K^{n}$ such that $\tuple{\cv, \delta \cv} \in W$.
\end{customenum}

\begin{definition}
We denote by 
\[
\Tdgp := \Td \cup \text{\ref{deep}}, \qquad
\Tdgpp := \Td \cup \text{\ref{wide}}.
\]
\end{definition}
We will show that both $\Tdgp$ and $\Tdgpp$ are axiomatizations for the model
completion of $\Td$.
Notice that the axiom scheme \ref{wide} deals with many variables at the same time,
but has only one iteration of the map $\delta$, while \ref{deep} deals with only one
variable at the same time, but many iteration of~$\delta$.

\begin{thm}\label{thm:MC} 
Assume that the theory T is model complete. 
Then the model completion $\Tdg$ of $\Td$ exists, and the theories
$\Tdgp$ and $\Tdgpp$ are two possible axiomatizations of $\Tdg$.
\end{thm}

\subsection{Proof preliminaries} In order to prove the main result we first introduce the following notation: given a polynomial $p(\x,y)$ we write
\[
p(\av, b) =^{y} 0 \iff\  
p(\av, b) = 0 \wedge \frac{\partial p}{\partial y} (\av, b) \neq 0.
\]

Let $S$ be a field of characteristic 0 and $\eps$ be a derivation on it. 
Let $I$ be a set of indexes (possibly, infinite).
Denote $\bar x \coloneqq \tuple{x_{i}: i \in I}$,
and $\bar y \coloneqq \tuple{y_{i}: i \in I}$.

\begin{definition}\label{def:delta-poly}
There exists a unique derivation
$S[\x] \to S[\x,\y]$, $p \mapsto p^{[\eps]}$ such that:
\[
\forall a \in S\quad a^{[\eps]} = \eps a, 
\forall  i \in I\quad x_{i}^{[\eps]} = y_{i};
\]
such derivation extends uniquely to a derivation
$S(\x) \to S(\x)[\y]$, $q \mapsto q^{[\eps]}$.

Moreover, the map $S(\x) \to S(\x)$ defined by
\[
q^{\eps} \coloneqq q^{[\eps]}(\x, 0)
\]
is the unique derivation on $S(\x)$ such that:
\[
\forall a \in S \quad a^{\eps} = \eps a, \qquad
\forall i \in I \quad x_{i}^{\eps} = 0;
\]
when $p \in S[\x]$, $p^{\eps}$ is the polynomial obtained by $p$ by applying
$\eps$ to each coefficient.
\end{definition}

\begin{remark}
For every  $q \in S(\x)$ and  $\av \in S^{n}$,
\begin{align}
q^{[\eps]}(\x, \y) &= q^{\eps}(\x) + \sum_{i \in I} \frac{\partial q}{\partial x_{i}}(\x) y_{i};\\
\eps(q(\av)) &= q^{[\eps]}(\av, \eps \av).
\end{align}
If moreover the field $S'$ is a field containing $S$ and $\eps': S(\x) \to S$ is a
derivation extending~$\eps$, then
\begin{equation}
\label{eq:1}
\eps'(q) = q^{[\eps]}(\x, \eps'(\x)) = q^{\eps}(\x) + \sum_{i \in I} \frac{\partial q}{\partial x_{i}}(\x) \eps'(x_{i}).
\end{equation}
\end{remark}

We will often also use the following fundamental fact, without further mentions
\begin{fact}\label{fact:extend-derivation}
Let $S'$ be a field  containing~$S$ (as in Definition~\ref{def:delta-poly}).
Let $\av \coloneqq \tuple{a_i: i \in I}$ be a (possibly, infinite) tuple of elements of $S'$ which are
algebraically independent over~$S$, and $\bv \coloneqq \tuple{b_{i}: i \in I}$ be a tuple of elements of~$S'$ (of
the same length as~$\av$).
Then, there exists a derivation $\eps'$ on $S'$ extending $\eps$ and such that
$\eps'(\av) = \bv$.
If moreover $\av$ is a transcendence basis of $S'$ over~$S$, then $\eps'$ is unique.
\end{fact}
\begin{proof}
\Wlog, $\av$ is a transcendence basis of $S'$ over~$S$.
By \cite[Ch.~II, \S17, Thm.~39]{ZS1}, 
there exists a unique derivation $\eps'': S(\av) \to S'$ extending $\eps$ and such that $\eps''(\av) = \bv$;
we can also prove it directly, by defining, for every $q \in S(\x)$,
\[
\eps''(q(\av)) \coloneqq q^{[\eps]}(\av, \bv).
\]

Given $c \in S'$, let $p(y) \in S(\av)[y]$ be the (monic) minimal polynomial of $c$
over~$S$.
Let $S'' \coloneqq S(c) \subseteq S'$.
Let
\begin{equation}
\label{eq:2}
d \coloneqq - p^{\eps''}(c)/ p'(c) \in S'.
\end{equation}
Any derivation on $\eps'$ on $S'$ extending $\eps''$ must satisfy $\eps'(c) = d$, and by
 \cite[Ch.~II, \S17, Thm.~39]{ZS1} again, 
there exists a unique derivation
$\eps''': S'' \to S'$ extending $\eps''$ and such that $\eps'''(c) = d$.

By iterating the above construction, we find a unique derivation $\eps'$ on $S'$
extending~$\eps''$.
\end{proof}

We need the following preliminary lemmas.
\begin{lemma}\label{polinomio} Let $\alpha(x, \y)$ be an L-formula and  
$(B, \delta) \models T^{\delta}$.  
Then there exists a  function $ \alpha^{[\eta]}$ definable in $T$  such that 
$\delta a = \alpha^{[\eta]}(a, \overline b, \delta \overline b)$, 
for every $a, \overline b \in B$ with  $B \models \alpha(a, \overline b)$ and $| \alpha( a, B) | <  \infty .$ 
\end{lemma}

\begin{proof}
Let $\alpha(x, \y)$ be an L-formula. Since $\K$ is algebraically bounded over $F$ and of
 characteristic 0,  there exist polynomials $p_1(x, \y), \ldots, p_k(x, \y) \in F[x,\y]$
 associated to the formula $\alpha(x, \y)$ and formulas 
$\beta_i(x, \y) ="p_i(x, \y) =^{x} 0"$
such that $T \vdash (\alpha(x, \y) ) \wedge |\alpha(x, \cdot )| < \infty ) \rightarrow \bigvee_{i=1}^{k} \beta_i(x, \y)$. 
Now we
can associate to each polynomial $p_i$ the partial function 
\[ f_i(x, \y, \delta \y) := \frac{\frac{\partial p_i}{\partial \y}\cdot \delta \y + p^{\eta}}{\frac{\partial p_i}{\partial x}},
\]
where $p^{\eta}$ is the polynomial defined in \ref{def:delta-poly} 
obtained by $p$ by applying $\eta$ to each
coefficients.

So now we have a total T-definable function $f(x, \y, \delta \y)$ whose graph is defined in the following way:
\begin{multline*} 
z = f(x, \y, \delta y ) \ \Leftrightarrow\   \Pa{ \beta_1(x, \y) \wedge  z=f_1(x, \y, \delta y)} \vee \Pa{\neg \beta_1(x, \y) \wedge
\beta_2(x, \y) \wedge  z=\\
=f_2(x, \y, \delta y)} \vee
\\ \vel \ldots \vee \Pa{\neg \beta_1(x, \y) \wedge \ldots \wedge \neg \beta_{k-1}(x, \y)  \wedge \beta_{k}(x, \y)  \wedge z=f_k(x,
  \y, \delta y)} \vee 
\\ \vel \Pa{\neg \beta_1(x, \y) \wedge \ldots \wedge \neg \beta_{k-1}(x, \y)  \wedge \neg \beta_{k}(x, \y)  \wedge  z= 0}.
\qedhere
\end{multline*}
\end{proof}

\begin{corollary}\label{derivata} For any T-definable function $f(\x)$ there exists a T-definable function $f^{[\eta]}$ such that $\delta(f(\x)) = f^{[\eta]}(\x, \Jet( \x)).$
\end{corollary}

\begin{lemma}\label{funzioni}
Le $t(\x)$ be an $L^{\delta}$-term. 
Then there is a T-definable function $f(\x, \y)$ such that 
$t(\x) = f(\x, \Jet(\x))$.
\end{lemma}

\begin{proof}
We prove by induction on the complexity of the term $t(\x).$
If $t(\x)$ is a variable it is trivial. Suppose that $t(\x) = h(s(\x)).$ By
induction there exists a $T$-definable function $g$ such that $s(\x) =g(
\x, \Jet(\x))$. 
If the function $h$ is in $L$ we can conclude. Otherwise $h = \delta$ and we obtain  $t(\x) = \delta(g(\x, \Jet(\x)))$. 
By Corollary  \ref{derivata} we conclude the proof. 
\end{proof} 

\begin{lemma}\label{senzaquantificatori}
Let $\phi(\x)$ be a quantifier free $L^{\delta}$-formula. Then there exists an L-formula 
$\psi$ such that $T^{\delta} \vdash \phi(\x) \leftrightarrow \psi(\x, \Jet(\x)).$
\end{lemma}

\begin{proof}
Follows from Lemma \ref{funzioni}.
\end{proof}

\subsection{Proof of Theorem \ref{thm:MC}}\label{sec:MC1-proof}

We can finally prove that both $\Tdgp$ and $\Tdgpp$ axiomatize $\Tdg$.
The proof is in three steps: firstly we show that $\Tdgpp \vdash \Tdgp$, and later we prove that the conditions \ref{en:Blum}
of Proposition~\ref{prop:MC} hold for $U = \Td$ and, more precisely, (a) holds for $U^{*}$ equal
to $\Tdgpp$ (\ie, that every model of
$\Td$ can be embedded in a model of $\Tdgpp$), and 
(b) for $U^{*}$ equal
 $\Tdgp$ (\ie, if $B \models \Td$ and $C \models \Tdgp$ have a common
substructure $A$, then every quantifier-free $\Ldelta(A)$-formula with one free variable having a
solution in $B$ also has a solution in $C$).

\begin{lemma}
$\Tdgpp \vdash \Tdgp$.
\end{lemma}
\begin{proof}
Let $Z \subseteq \K^{n+1}$ be $L(\K)$-definable such that $\Pi_{n}(Z)$ is large.
Define
\[
W := \set{\tuple{\x, \y} \in \K^{n} \times \K^{n}:
\tuple{\x, y_{n}} \in Z \wedge \bigwedge_{i=1}^{n-1} y_{i} = x_{i+1}}.
\] 
Clearly, $\Pi_{n}(W) = \Pi_{n}(Z)$, and therefore $\Pi_{n}(W)$ is large.
By \ref{wide}, there exists  $\cv \in \K^{n}$ such that
$\tuple{\cv, \delta \cv} \in W$.
Then, $\Jetdn(c_{1}) \in Z$.
\end{proof}

\begin{lemma}\label{lem:embedding}
Let $(A, \delta) \models \Td$.
Let $Z \subseteq A^{n} \times A^{n}$ be $L$-definable with parameters in $A$, such that
$\Pi_{n}(Z)$ is large.
Then, there exists $\tuple{B, \eps} \supseteq \tuple{A, \delta}$ and $\bv \in B^{n}$ such that
$B \succeq A$,
$\tuple{B,\eps} \models \Td$, and 
$\tuple{\bv, \eps \bv} \in Z_{B}$ (the interpretation of $Z$ in B).
\end{lemma}
\begin{proof}
Let $B \succ A$ (as $L$-structures) be such that $B$ is $\card{A}^{+}$-saturated.
By definition of dimension, there exists $\bv \in \Pi_{n}(Z_{B})$ which is
algebraically independent over~$A$.
Let $\dv \in B^{n}$ be such that $\tuple{\bv, \delta \bv} \in Z_{B}$.
Let $\eps$ be any derivation on $B$ which extends $\delta$ and such that, by Fact \ref{fact:extend-derivation}, $\eps (\bv) = \dv.$
\end{proof}

\begin{lemma} Let $\tuple{B,\delta} \models \Td$, $\tuple{C,\delta} \models \Tdgp$,
and $\tuple{A,\delta}$ be an $L^{\delta}$-substructures of both models,
such that $B$ and $C$ have the same $L(A)$-theory.
Let  $b \in B$ be such that $\tuple{B, \delta} \models \theta(b)$,
where $\theta(x)$ is a quantifier free $L^{\delta}$-formula with parameters in~$A$. 
Then, there exists $c \in C$ such that $\tuple{C,\delta} \models \theta(c)$.
\end{lemma}

\begin{proof}
By Lemma \ref{senzaquantificatori}  there exist $n \in \N$ and an $L(A)$-formula $\psi$  such
that $\theta(x) = \psi (\Jetdn(x))$. 

Let $Y^{B} := \psi(B) = \set{\dv \in B^{n+1}: B \models \psi(\dv)}$, and $Y^{C} := \psi(C)$. 

Let $d$ be the smallest integer such that $\delta^d (b)$ is algebraically dependent
over $\Jet^{d-1}(b) \cup A$ (or $d = + \infty$ if $\Jetdinf(b)$ is algebraically
independent over~$A$).
We distinguish two cases:

1) $d \geq n$: in this case, $\Pi_{n}(Y^C)$ is large 
because $\Jet^{n-1}(b) \in \Pi_{n}(Y^{B})$, therefore,
by \ref{deep}, there exists $c \in C$ such that $C \models \theta(\Jetdn(c)).$

2) $d < n$: this means that $\delta^{d}b \in \acl(\Jet^{d-1}(b)),$ so there exists 
polynomial $p(\y, x) \in A[\y,x]$ such that $p(\Jet^{d-1}(b), \delta^{d} b) =^{x} 0$.
By Lemma \ref{polinomio} there exist  $L(A)$-definable functions 
$f_{d+1}, f_{d+2}, \ldots, f_n$  such that $\delta^{i} = f_i(\Jet^{d}(b))$ where 
$i = d + 1, d + 2, \ldots, n.$ 
Let \[
Z^{B} := \{\y \in B^{d+1}  : p(\y) = ^{y_{d+1}} 0  \wedge \theta(\y, f_{d+1}(\y), \ldots,
f_{n}(\y))\}.
\]
Notice that $\Pi_{d}(Z^{C})$ is large, because  $\Jet^{d-1}(b) \in \Pi_{d}(Z^{B})$, 
and therefore by axiom \ref{deep} there exists $c \in C$ such that $\Jet^{d}(c) \in Z^{C}$ and so $\Jet^{n}(c) \in Y^{C}$.
\end{proof}

\subsection{Corollaries}

\begin{corollary}
Assume that $T$ eliminates quantifiers.
Then,  $\Tdgp$ and $\Tdgpp$ are axiomatizations for the model completion $\Tdg$
of $\Td$.

Moreover, $\Tdg$ admits elimination of quantifiers, and for every $\Ldelta$-formula
$\alpha(\x)$ there exists a quantifier-free $L$-formula $\beta(\y)$ such that
\[
\Tdg \models \forall \x \, \Pa{\alpha(\x) \leftrightarrow \beta(\Jet  (\x) )}. 
\]

Finally, $\Tdg$ is complete.
\end{corollary}

\begin{corollary}
Assume that $T$ is model complete.
Then,  $\Tdgp$ and $\Tdgpp$ are axiomatizations for the model completion $\Tdg$
of $\Td$.
\end{corollary}

The next corollary is without any further assumptions on $T$.
\begin{corollary}
$\Tdgp$ and $\Tdgpp$ are equivalent consistent theories (which we denote by
$\Tdg$).

Moreover, for every $\Ldelta$-formula
$\alpha(\x)$ there exists an $L$-formula $\beta(\y)$ such that
\[
\Tdg \models \forall \x \, \Pa{\alpha(\x) \leftrightarrow \beta(\Jet \x)}. 
\]

Finally, $\Tdg$ is complete.
\end{corollary}

\section{Several non-commuting derivations}\label{sec:nc}

We analyze first the case when there are several not commuting derivations
$\delta_{1}, \dotsc, \delta_{k}$ because it is simpler in terms of axiomatization, as we observed in Section \ref{storico}, and later  in Section \ref{sec:commute} we examine the harder case of commuting derivations.

Let $\deltabar := \tuple{\delta_{1}, \dotsc, \delta_{k}}$.
Let $\eta_{1}, \dotsc, \eta_{k}$ be  derivations on~$F$.
We denote by $\Tdn$ the $\Ldb$-expansion of $T$ saying that each $\delta_{i}$ is a
derivation and that $\delta_{i}$ extends $\eta_{i}$ for $i \leq k$.

\begin{thm}
Assume that $T$ is model complete.
Then, $\Tdn$ has a model completion $\Tdng$. 
\end{thm}

To give the axioms for $\Tdng$ we need some more definitions and notations.
We fix $\tuple{\K, \deltabar} \models \Tdn$.

Let $\Gamma$ be the free non commutative monoid generated by $\deltabar$, with the
canonical partial order~$\preceq$ given by $\beta \preceq \alpha\beta$, for all $\alpha, \beta \in \Gamma.$
We fix the total order on $\Gamma$, given by
\[
\theta \leq \theta'  \Leftrightarrow  \abs \theta < \abs {\theta'} \vel
\Pa {\abs \theta = \abs {\theta'} \et \theta <_{lex} \theta'},
\]
where $<_{lex}$ is the lexicographic order, and
$\abs{\theta}$ is the length of $\theta$ as a word in the alphabet~$\deltabar$.

\begin{remark}
$\preceq$ is a well-founded partial order on $\Gamma$, but it is not a well-partial-order
(\ie, there exist infinite anti-chains).
\end{remark}

\begin{remark}
\begin{enumerate}
\item As an ordered set, $\tuple{\Gamma, \leq}$ is isomorphic to $\tuple{\N, \leq}$;
\item $\emptyset$ (\ie, the empty word, corresponding to the identity function on $\K$)
 is the minimum of $\Gamma$;
\item If $\alpha \preceq \beta$, then $\alpha \leq \beta$;
\item If $\alpha \leq \beta$, then $\gamma \alpha \leq \gamma \beta$ and
$\alpha \gamma \leq \beta \gamma$, for all $\gamma \in \Gamma$.
\end{enumerate}
\end{remark}

For every variable $x$ and every $\gamma \in \Gamma$ we introduce the variable $x_{\gamma}$.
Given $V \subseteq \Gamma$, we denote $x_{V} \coloneqq \tuple{x_{\gamma}: \gamma \in V}$
 and $a^{V} \coloneqq \tuple{\gamma a: \gamma \in V}$. We remark that $a^V$ is an analogue of the notion of Jet in one derivation, i.e. $\Jet^n(a) = a^{\{0, 1, \ldots, n\}}.$
Moreover, we denote $\Pi_{A}$ the projection from $\K^{B}$ to $\K^{A}$
 (for some $A, B \subseteq \Gamma$ and $B \supseteq A$), mapping $\tuple{a_{\mu}: \mu \in B}$ to
$\tuple{a_{\mu}: \mu \in A}$.

We give now two alternative axiomatizations for $\Tdng$.

\begin{enumerate}[label=(nc-\ctext*), 
  ref=(\texttt{nc-\ctext*}),
  font= \texttt]
\item\label{nc-deep}
Let $\tildefree \subset \Gamma$ be finite and $\preceq$-initial.
Let $\thepivot \subseteq \tildefree$ be the set of $\preceq$-maximal elements of $\tildefree$,
and $\thefree \coloneqq \tildefree \setminus \thepivot$.
Let $Z \subseteq \K^{\tildefree}$ be $L(A)$-definable.
If $\Pi_{\thefree }(Z)$ is large, then there exists $c \in \K$ such that
$c^{\tildefree} \in Z$.
\item\label{nc-wide}
Let  $W \subseteq \K^{n} \times \K^{k \times n} $ be $L(\K)$-definable, such that $\Pi_{n}(W)$ is
large.
Then, there exists 
$\cv \in \K^{n}$ such that $\tuple{\cv, \delta_{1} \cv, \dotsc, \delta_{k} \cv} \in W$.
\end{enumerate}

\begin{definition}
We denote by 
\[
\Tdngp := \Td \cup \text{\ref{nc-deep}}, \qquad
\Tdngpp := \Td \cup \text{\ref{nc-wide}}.
\]
\end{definition}

\begin{thm}\label{thm:ncMC}
\begin{enumerate}
\item $\Tdngp$ and $\Tdngpp$ are consistent and equivalent to each other.
\item If $T$ is model-complete, then 
the model completion $\Tdng$ of $\Tdn$ exists, and the theories
$\Tdngp$ and $\Tdngpp$ are two possible axiomatizations of $\Tdng$.
\item If $T$ eliminates quantifiers, then $\Tdng$ eliminates quantifiers.
\item For every $\Ldeltabar$-formula $\alpha(\x)$ there exists an $L$-formula $\beta(\x)$
such that
\[
\Tdng \models \forall \x\ \Pa{\alpha(\x) \leftrightarrow \beta(\x^{\Gamma})}.
\]
\end{enumerate}
\end{thm}

For the proof, we proceed as in \S\ref{sec:MC1-proof}, i.e. it is in three steps:
\begin{lemma}
$\Tdngpp \vdash \Tdngp$.
\end{lemma}
\begin{proof}
Let $Z, \thefree, \thepivot, \tildefree$ be as in \ref{nc-deep}.
\begin{claim}
\Wlog, we may assume that $\thepivot$ is equal to the set of $\preceq$-minimal
elements of $\Gamma \setminus \thefree$.
\end{claim}
In fact, let  $\thepivot'$ be  the set of $\preceq$-minimal
elements of $\Gamma \setminus \thefree$; notice that $\thepivot \subseteq \thepivot'$. We can replace $\tildefree$ with $\tildefree' \coloneqq \tildefree \cup \thepivot'$,
and $Z$ with $Z' \coloneqq \Pi^{-1}(Z)$, where the function $\Pi$ is defined as:
\functiondef%
{\Pi}%
{\K^{\tildefree'}}%
{\K^{\tildefree}}%
{\x}%
{\tuple{x_{\mu}: \mu \in \tildefree}.}

Then, $\Pi_{\thefree}(Z') = \Pi_{\thefree}(Z)$, and if $a^{\tildefree'} \in Z'$,
then $a^{\tildefree} \in Z$.

\smallskip

We introduce variables $x_{0}, x_{1}, \dotsc, x_{k}$ and corresponding
variable $x_{i,\gamma}$, 
which for readability we denote by $x(i,\gamma)$  such that $0 \leq i \leq k, \gamma \in \Gamma$.
For brevity, we denote 
\[ \x \coloneqq \tuple{x(i,\gamma): 0 \leq i \leq k, \gamma \in \tildefree} \qquad \text{and} \qquad
\x_{i} \coloneqq {\tuple{x(i,\gamma): \gamma \in \tildefree}}, \quad i= 0, \dotsc, k.
\]

We also denote
\functiondef%
{\Pi_{0}}%
{(\K^{\tildefree})^{k+1}}%
{\K^{\tildefree}}%
{\x} 
{\x_0} 
For each $\pi \in \thepivot$, we choose 
$\mu_{\pi} \in \thefree$ and $i_{\pi} \in \set{1,   \dotsc, k}$ such that $\delta_{i_{\pi}} \mu_{i} = \pi$.
Moreover, given $\av \in (\K^{\thefree})^{k+1}$,
we define $\av' \in K^{\tildefree}$ as the tuple with coordinates
\[a'_{\gamma} \coloneqq 
\begin{cases}
a(0, \gamma) & \text{ if } \gamma \in \thefree\\
a(i_{\gamma}, \mu_{\gamma})
& \text{  if } \gamma \in \thepivot.
\end{cases}
\]
We define
\[
W \coloneqq \set{\tuple{\av \in (\K^{\thefree})^{k+1}}:
\av' \in Z 
\et a(i, \gamma) = a(0, \delta_{i}\gamma) \text{ for } i =1, \dotsc,k
\text{ and } \gamma \in \thefree
}.
\]
Notice that $\Pi_{0}(W)$ is equal to $\Pi_{\thefree}(Z)$, and therefore it is large.
Thus, by \ref{nc-wide}, there exists $\av \in \K^{\thefree}$ such that
$\tuple{\av, \delta_{1}(\av), \dotsc, \delta_{k}(\av)} \in W$.
Finally, taking $a \coloneqq a(0,\emptyset)$, we get
$a^{\tildefree} \in Z$.
\end{proof}

\begin{lemma}
Let $(A, \deltab) \models \Tdn$.
Let $Z \subseteq A^{n} \times (A^{n})^{k}$ be $L$-definable with parameters in $A$, such that
$\Pi_{n}(Z)$ is large.
Then, there exists $\tuple{B, \epsb} \supseteq \tuple{A, \deltab}$ and $\bv \in B^{n}$ such that
$B \succeq A$,
$\tuple{B,\epsb} \models \Tdn$, and 
$\tuple{\bv, \epsb \bv} \in Z_{B}$.
\end{lemma}
\begin{proof}
Same proof as for Lemma~\ref{lem:embedding}.
\end{proof}

\begin{lemma} Let $\tuple{B,\deltab} \models \Tdn$, $\tuple{C,\deltab} \models \Tdngp$,
and $\tuple{A,\deltab}$ be an $L(\deltab)$\hyph{}substructures of both models,
such that $B$ and $C$ have the same $L(A)$-theory.
Let  $b \in B$ be such that $\tuple{B, \deltab} \models \theta(b)$,
where $\theta(x)$ is a quantifier free $L(\deltab)$-formula with parameters in~$A$. 
Then, there exists $c \in C$ such that $\tuple{C, \deltab} \models \theta(c)$.
\end{lemma}
\begin{proof}
By Lemma \ref{senzaquantificatori}  there exists $U$ finite subset of $\Gamma$ and
an $L(A)$-formula $\psi(\y)$  such that $U$ is $\preceq$-initial and
that $\Tdn \models \theta(x) = \psi (x^{U})$. 
Let $Y^{B} \coloneqq \Psi(B)$ and $Y^{C} \coloneqq \psi(C)$.
Let 
\[
\thefree \coloneqq \set{\gamma \in U: \gamma b \notin \acl(A, b^{U < \gamma})}, \mbox{ where we denote
} b^{U < \gamma}  \coloneqq \tuple{\mu b: \mu < \gamma \wedge \mu \in U}.
\]
Define $\thebounded \coloneqq \Gamma \setminus \thefree$ and
$\thepivot$ be the set of $\preceq$-minimal elements of $\thebounded$
(notice that $\thepivot$ might be infinite).
As usual, define $\tildefree \coloneqq \thefree \cup\ \thepivot$.

For every $\gamma \in \Gamma$ there exists $q_{\gamma} \in A(x_{\tildefree \leq \gamma})$ such that
$\gamma b = q_{\gamma}(b^{\tildefree \leq \gamma})$. 
Let $\beta$ be the following $L(A)$-formula:
\[
\beta(x_{\tildefree}) \coloneqq \psi(q_{\gamma}(x_{\tildefree}): \gamma \in U).
\]
Notice that $\tuple{B, \deltab} \models \beta(b^{\tildefree})$.
Let $\tildefree_{0} \subseteq \tildefree$ be the set of indexes of the variables of~$\beta$: 
\wloG, we may assume that
$\tildefree_{0}$ is a $\preceq$-initial subset of~$\Gamma$.
Let $\thepivot_{0}$ be the set of $\preceq$-maximal elements of $\tildefree_{0}$.
Define
\[
Z \coloneqq \set{\dv \in \K^{\tildefree_{0}}: \tuple{B, \deltab} \models \beta(\dv)}.
\]
Notice that $\Pi_{\thefree_{0}}(Z)$ contains $b^{\thefree_{0}}$,
and therefore it is large.
Thus, by~\ref{nc-deep}, there exists $c \in C$ such that
$c^{\tildefree_{0}} \in Z$, and therefore $c^{U}$ satisfies~$\psi$.
\end{proof}


\section{Several commuting derivations}\label{sec:commute}
We now deal with the case when there are several \textbf{commuting} derivations
$\delta_{1}, \dotsc, \delta_{k}$. The technique used here for the treatment of the study of several derivations are a variant of \cite{FK}. In particular, we avoid as much as possible the algebraic approach in \cite{Kolchin} based on autoreduced sets.

Let $\deltabar := \tuple{\delta_{1}, \dotsc, \delta_{k}}$ and 
let $\eta_{1}, \dotsc, \eta_{k}$ be commuting derivations on~$F$.
Let $\Tdb$ be the $\Ldb$-expansion of $T$ saying that each $\delta_{i}$ is a
derivation, that $\delta_{i}$ extends $\eta_{i}$ for $i \leq k$,
and that $\delta_{i} \circ \delta_{j} = \delta_{j} \circ \delta_{i}$, for $i, j \leq k$.

\begin{thm}
Assume that $T$ is model complete.
Then, $\Tdb$ has a model completion $\Tdbg$. 
\end{thm}

\subsection{Configurations}
In order to give axioms for~$\Tdbg$, 
we first need some more definitions and notations.
We fix $\tuple{\K, \deltabar} \models \Tdb$.
We denote by $K$ the field underlying~$\K$.

Let $\Theta$ be the free commutative monoid generated by $\deltabar$, with the
canonical partial order~$\preceq$ (notice that $\Theta$ is isomorphic to~$\N^{k}$).

Notice that $\Theta$ is, canonically, a quotient of the free monoid $\Gamma$. For every $v
\in \Gamma$ we denote by $[v]\in \Theta$  the equivalence class of~$v$.

We fix the total order on $\Theta$, given by
\[
\theta \leq \theta' \mbox{ iff } \abs \theta < \abs {\theta'} \vel
\Pa {\abs \theta = \abs {\theta'} \et \theta <_{lex} \theta'},
\]
where $<_{lex}$ is the lexicographic order, and
$\abs{\tuple{\delta_{1}^{n_{1}} \cdots \delta_{k}^{n_{k}}}} \coloneqq n_{1} + \dots + n_{k}$.

Given $a \in \K$ and $\theta \in \Theta$, we denote by 
$a^{< \theta} \coloneqq \tuple{\mu a: \mu < \theta}$, 
and similarly $a^{\leq \theta} \coloneqq \tuple{\mu a: \mu \leq \theta}$, 
and $a^{\Theta} \coloneqq \tuple{\mu a: \mu \in \Theta}$.
Moreover, for each $\theta \in \Theta$ we have a variable $x_{\theta}$, and we denote
$x_{<\theta} \coloneqq \tuple{x_{\mu}:  \mu < \theta}$.
Moreover, given a set $A \subseteq \Theta$, we denote $x_{A} \coloneqq \tuple{x_{\theta}: \theta \in A}$,
and $x_{A \leq \mu} \coloneqq \tuple{x_{\nu}: \nu \in A \et \nu \leq \mu}$.
Given a rational function $q \in K(x_{\Theta})$ and $\av \in K^{\Theta}$, we denote by
\[
\frac{\partial q}{\partial \mu} \coloneqq \frac{\partial q}{\partial x_{\mu}} \qquad \text{and} \qquad
q(\av) =^{\mu} 0 \mbox{ iff } q(\av) = 0 \et \frac{\partial q}{\partial \mu}(\av) \neq 0.
\]

\medskip

Let $K_{0}$ be a differential subfield of $K$ (\ie, such that $\deltabar(K_{0})
\subseteq K_{0}$).

A \intro{\posture}  $\condition$
with parameters in~$K_{0}$ is given by the following data.

1) A $\preceq$-anti-chain  $\thepivot \subset \Theta$. Notice that, by Dickson's Lemma,
$\thepivot$ must be finite.

We distinguish two sets:
\begin{itemize}
\item $\thebounded \coloneqq \set{\mu \in \Theta: \exists \pi \in \thepivot\, \mu \succeq \pi}$, the set of leaders;
\footnote{We take the terminology from \cite{Pierce14}, who in turns borrows it
  from Ritt. Then $\thepivot$ is the set of minimal leaders.}
\item $\thefree \coloneqq \Theta \setminus \thebounded = \set{\mu \in \Theta: \forall \pi \in \thepivot\, \mu \nsucceq \pi}$, the
set of free elements.\\
Moreover, we define:
\item $\tildefree \coloneqq \thefree \cup \thepivot$. 
\end{itemize}
Notice that $\thepivot$ is the set of $\preceq$-minimal elements of~$\thebounded$.
We assume that $\thefree$ is non-empty (equivalently, that $0 \in \thefree$).

2) For every $\pi \in \thepivot$ we are given a nonzero polynomial
$p_{\pi} \in K_{0}[x_{\pi}, x_{\thefree < \pi}]$ which depends on~$x_{\pi}$
(\ie, its degree in $x_{\pi}$ is nonzero).

Consider the quasi-affine variety (defined over $K_{0}$)
\[
W_{0} \coloneqq \set{x_{\tildefree} \in K^{\tildefree}: \bigwedge_{\pi \in \thepivot}
  p_{\pi}(x_{\tildefree}) = ^{\pi} 0},
\]
(by quasi-affine variety we simply mean a subset of $K^{n}$ which is locally
closed in the Zariski topology. We don't consider its spectrum).

3) Finally, we are given $W \subseteq W_{0}$, which is Zariski closed in $W_{0}$, such
that $W$
is defined (as a quasi-affine variety) over~$K_{0}$
and such that
$\Pi_{\thefree}(W)$ is large (where $\Pi_{\thefree}: K^{\tildefree} \to
K^{\thefree}$ is the canonical projection).

\smallskip
For the remainder of this section, 
we are given a \posture:
\[
\condition :=  (W; p_{\pi}: \pi \in \thepivot).
\]

\medskip

We want to impose some commutativity on $\condition$.




We define now some induced data.

For every $\alpha \in \Theta$ we define a rational function 
$f_{\alpha} \in K_{0}(x_{\tildefree})$ and a tuple of rational functions~$F_{\alpha}$.
\begin{itemize}
\item If $\alpha \in \tildefree$, then $f_{\alpha} \coloneqq x_{\alpha}$ and $F_{\alpha} =
\set{f_{\alpha}}$.
\item 
For every $\pi \in \thepivot$ and $v \in \Gamma$, we define
$f_{v, \pi}$ by induction on $v$, and then 
\[
F_{\alpha} := \set{f_{v, \pi}: v \in \Gamma, \pi \in \thepivot, [v]\pi = \alpha},
\]
and $f_{\alpha}$ is an arbitrary function in~$F_{\alpha}$.\\
If $v = 0$, then $f_{v,\pi} \coloneqq f_{\pi} = x_{\pi}$.\\
If $v = \delta$ for some $\delta \in \deltabar$, define
\begin{equation}
\label{eq:3}
f_{\delta, \pi} \coloneqq - \frac{p_{\pi}^{\delta} + \sum_{\mu \in \thefree} \frac{\partial p_{\pi}}{\partial \mu}\cdot f_{\delta
  \mu}}  {\frac{\partial p_{\pi}}{\partial \pi}}.
\end{equation}
If $v = \delta w$ with $0 \neq w \in \Gamma$, define
\[
f_{\delta w, \pi} \coloneqq f_{w, \pi}^{\delta} + 
 \sum_{\mu \in \tildefree} \frac{\partial f_{w,\pi}}{ \partial \mu} \cdot f_{ \delta, \mu},
\]
where $f_{\delta, \mu} \coloneqq f_{\delta \mu}$ when $\mu \in \thefree$.
\end{itemize}

We also define $g_{\alpha}$ and $g_{w, \pi}$ as the restriction to $W$ of
$f_{\alpha}$ and $f_{w, \pi}$, respectively, and
$G_{\alpha}$ the family $\{ g_{w, \pi}: w \in \Gamma, \pi \in \thepivot, [w]\pi = \alpha \}$.

Let $\theta\in \Theta$ be the $\preceq$-l.u.b.\ of $\thepivot$
and $d$ be the dimension of~$W$.
We say that $\condition$ \intro{commutes} at $\alpha \in \Theta$ if, for every $g,g' \in G_{\alpha}$,
$g$ and $g'$ coincide outside a subset of $W$
of dimension less than $d$
(\ie, they coincide ``almost everywhere'' on~$W$).
We say that $\condition$ commutes locally if it commutes at every $\alpha \leq \theta$, and
it commutes globally if it commutes at every $\alpha \in \Theta$.

The main result that makes the machinery work is the following.
\begin{thm}\label{thm:local-global-commute}
$\condition$ commutes locally iff it commutes globally.
\end{thm}
We will say then that $\condition$ commutes if it commutes locally
(equivalently, globally).

In order to prove the above theorem, we need some preliminary definitions and
results.  
It suffices to prove the theorem for every $K_{0}$-irreducible component of~$W$ 
of dimension~$d$.
Thus, \wloG we may \textbf{assume} for the remainder of this subsection that $W$
is $K_{0}$-irreducible.
Then, the condition that $W$ commutes at a certain $\alpha$ becomes that $G_{\alpha}$ is a
singleton.

We denote by $K_{0}[W]$ the ring of regular functions on $W$ (that is, the
restriction to $W$ of polynomial maps $K_{0}^{\tildefree} \to K_{0}$).
Since we are assuming that $W$ is $K_{0}$-irreducible, 
$K_{0}[W]$ is an integral domain.
Therefore, we can consider its
fraction field~$K_{0}(W)$.
An element of $K_{0}(W)$ is the restriction to $W$ of a rational function in
$K_{0}(x_{\tildefree})$; in particular, the functions $g_{\alpha}$ and $g_{w, \pi}$ are 
in~$K_{0}(W)$.

Given $\delta \in \deltabar$, we define the following function:
\begin{functiondef}
{R^{\delta}_{0}}
{K_0(x_{\tildefree})} 
{K_0(x_{\tildefree})}
h 
{h^{\delta} + \sum_{\mu \in \tildefree} \frac{\partial h}{\partial \mu} \cdot f_{\delta, \mu}.}
\end{functiondef}
Notice that 
$R^{\delta}_{0}$ is the unique derivation extending $\delta$  such that 
$R^{\delta}_{0}(x_{\mu}) = f_{\delta, \mu}$, for every $\mu \in \tildefree$.
Given $w = w_{1} \cdots w_{\ell} \in \Gamma$, we can define 
$R^{w}_{0}: K_{0}(x_{\tildefree}) \to K_{0}(x_{\tildefree})$ as the
composition $R_{0}^{w} \coloneqq R_{0}^{w_{1}} \circ \dots \circ R_{0}^{w_{\ell}}$.

We have, for every $v,w \in \Gamma$ and every $\pi \in \thepivot$,
\begin{equation}
\label{eq:8}
R^{w}_{0}(f_{v, \pi}) = f_{wv, \pi}.
\end{equation}

If we restrict $R^{\delta}_{0}$ to $K_{0}(x_{\thefree}) $ and compose with the
restriction to~$W$, we obtain a derivation
$R^{\delta}_{1}: K_{0}(x_{\thefree}) \to K_{0}(W)$. 
An equivalent definition is that
$R^{\delta}_{1}$ is the unique derivation extending $\delta$ such that
$R^{\delta}_{1}(x_{\mu}) = g_{\delta \mu}$ for every $\mu \in \thefree$.
Finally, observe that $K_{0}(W)$ is an algebraic extension of 
$K_{0}(x_{\thefree})$.
So, $R^{\delta}_{1}$ extends uniquely to a derivation $R^{\delta}$ from
$K_{0}(W)$ to the algebraic closure of~$K_{0}(W)$.

We will consider the objects $x_{\mu}$ both as variables and as functions.
Observe that, for every $\mu \in \tildefree$, 
$g_{\mu}$ is the restriction of $x_{\mu}$ to~$W$.

\begin{remark}\label{rem:Rdelta-rational}
Let $f \in K_{0}(x_{\tildefree})$ and $h \in K_{0}(W)$ be the restriction of $f$ to $W$.
Then, we have $h = f(g_{\tildefree})$ (we are seeing $f$ as a rational function).
Therefore,
\[
R^{\delta}(h) = f^{\delta}(g_{\tildefree}) + 
\sum_{\mu \in \tildefree} \frac{\partial f}{\partial \mu}\rest_{W} \cdot R^{\delta}(g_{\mu}) = 
f^{\delta} \rest_W + \sum_{\mu \in \tildefree} \frac{\partial f}{\partial \mu}\rest_{W} \cdot  g_{\delta, \mu}.
\]
\end{remark}

\begin{lemma}\label{lem:Rdelta-one}
Let $\mu \in \tildefree$.
Then, 
\begin{equation}
\label{eq:7}
R^{\delta}(g_{\mu}) = g_{\delta, \mu}.
\end{equation}
\end{lemma}
\begin{proof}
If $\mu \in \thefree$, the conclusion follows by definition of $R^{\delta}_{1}$.

If $\mu = \pi \in \thepivot$, observe that $g_{\pi}$ satisfies the algebraic condition
\[
p_{\pi}(g_{\pi}, g_{\thefree<\pi}) =^{\pi} 0.
\]
Therefore,
\[
R^{\delta}(g_{\pi}) = - \frac{p_{\pi}^{\delta} + \sum_{\beta \in \thefree < \pi}
\frac{\partial p_{\pi}}{\partial \beta} \cdot R^{\delta}(g_{\beta})}
{\frac{\partial p_{\pi}}{\partial \pi}}(g_{\pi}, g_{\thefree<\pi}) = 
f_{\delta, \pi} \rest_{W} =
g_{\delta, \pi}.  \qedhere
\]
\end{proof}

\begin{remark}
The image of $R^{\delta}$ is already in $K_{0}(W)$ (no need to take the algebraic
closure). Indeed, as a $K_{0}$-algebra, $K_{0}(W)$ is generated by
$(g_{\mu}: \mu \in \tildefree)$.
Thus, it suffices to show that $R^{\delta}(g_{\mu}) \in K_{0}(W)$, which
follows from Lemma~\ref{lem:Rdelta-one}.
\end{remark}
\begin{lemma}\label{lem:Rw}
For all $w \in \Gamma,$ $R^w(g_{v, \pi}) = g_{wv, \pi}$.
\end{lemma}
\begin{proof}
Apply  \eqref{eq:7} and \eqref{eq:8}.
%
%
%
\end{proof}

\begin{lemma}
Let $\delta, \eps \in \deltabar$, $f \in K_{0}(x_{\tildefree})$, and  $h \in K_{0}(W)$ be the
restriction of $f$ to~$W$. Then,
\begin{equation}
\label{eq:5}
[R^{\eps}, R^{\delta}]h = \sum_{\mu \in \tildefree} 
\frac{\partial f }{\partial \mu}\rest _{W} \cdot (R^{\eps} g_{\delta, \mu} - R^{\delta}  g_{\eps, \mu}).
\end{equation}
\end{lemma}
\begin{proof}
\textsc{First proof)}
Let $\lambda = [\eps, \delta]$ be the Lie bracket of $\eps$ and $\delta$, and
$R^{\lambda} = [R^{\eps}, R^{\delta}]$ be corresponding Lie bracket.
We can write $h = f(g_{\tildefree})$, where we see $f$ as a rational function: therefore, since $R^{\lambda}$ is a
derivation, we have
\[
R^{\lambda} h = h^{\lambda} + 
\sum_{\mu \in \tildefree} \frac{\partial f}{\partial \mu}(g_{\tildefree}) \cdot R^{\lambda}(g_{\mu}) =
\sum_{\mu \in \tildefree} 
\frac{\partial f }{\partial \mu}\rest _{W} \cdot (R^{\eps} g_{\delta, \mu} - R^{\delta}  g_{\eps, \mu}).
\]

\textsc{Second proof)}
Since both the LHS and the RHS of \eqref{eq:5} define a derivation on $K_{0}(W)$, it
suffices to show that they are equal when $f = x_{\mu}$ for some $\mu \in \thefree$.

Then, RHS is equal to
$R^{\eps} g_{\delta, \mu} - R^{\delta} g_{\eps, \mu}$, which is equal to 
$R^{\eps} R^{\delta} g_{\mu} - R^{\delta} R^{\eps} g_{\mu}$, \ie the LHS.
\end{proof}

We can finally prove the theorem.
\begin{proof}[Proof of Thm.~\ref{thm:local-global-commute}]
Assume that $\condition$ commutes locally (and that $W$ is $K_{0}$-irreducible).
Let $\alpha \in \thebounded$:
we show, by induction on~$\alpha$, that  $\condition$  commutes
at~$\alpha$.
Let $\pi, \pi' \in \thepivot$ and $w, w' \in \Gamma$ be such that
$[w] \pi = [w'] \pi' = \alpha$: we need to show that $g_{w, \pi} = g_{w', \pi'}$.

If $\pi = \pi'$, we can reduce to the case when
$w= v \delta \eps u$ and $w' = v \eps \delta u$ for some $\delta, \eps \in \deltabar$,
$v, u \in \Gamma$.\\

If $u \neq 0$, we have, by inductive hypothesis, that
$g_{v \delta \eps, \pi} = g_{v \eps \delta,  \pi}$ and therefore
\[
g_{w, \pi} = R^{u} g_{v \delta \eps, \pi} = R^{u} g_{v \eps \delta,  \pi} = g_{w', \pi},
\]
which is the thesis.

If instead $u = 0$,
we have
\[
g_{w, \pi} - g_{w', \pi} = [R^{\eps}, R^{\delta}] g_{v, \pi} =
\sum_{\mu \in \tildefree} 
\frac{\partial f_{v, \pi} }{\partial \mu}\rest _{W} \cdot (R^{\eps} g_{\delta, \mu} - R^{\delta}  g_{\eps, \mu}).
\]
Fix some $\mu$ in the above sum. Notice that $\mu < [v]\delta$, and therefore
$\delta \eps \mu = \eps \delta \mu < \alpha$.
Therefore, by inductive hypothesis,
$R^{\eps} g_{\delta, \mu} = g_{\eps \delta \mu} = R^{\delta}  g_{\eps, \mu}$. Thus, all summands are 
$0$ and $g_{w, \pi} - g_{w', \pi} = 0$, and we are done.

If $\pi \neq \pi'$, let $\beta \coloneqq \pi \vel \pi'$; by definition, $\beta \leq \theta$, and therefore,
by assumption, $G_{\beta} = \set{g_{\beta}}$.
Let $u, u', w \in \Gamma$ be such that 
\[
\beta = [u] \pi = [u']\pi', \quad \alpha = [w] \beta.
\]
By the previous case, we have
\[
g_{v, \pi} = g_{w u, \pi} = R^{w} g_{u, \pi} = R^{w} g_{\beta} 
= R^{w} g_{u', \pi'} = g_{w u', \pi'} = 
g_{v' \pi'}. 
\qedhere
\]
\end{proof}

\begin{remark}\label{rem:R-comm}
If $\condition$ commutes globally, then the derivations $R^{\delta_{1}}, \dotsc,
R^{\delta_{k}}$ commute with each other.
\end{remark}

The following remark motivates the definition of the functions~$g_{\mu}$.
\begin{remark}\label{rem:g-iteration}
Let $b \in K$ be such that $b^{\tildefree} \in W$.
Then, for every $w \in \Gamma$ and $\pi \in \thepivot$
\[
b^{[w]\pi} = g_{w, \pi}(b).
\]
Therefore, for every $\mu \in \Theta$,
\[
b^{\mu} = g_{\mu}(b).
\]
\end{remark}
\begin{proof}
By induction on~$w$.
\end{proof}

\begin{corollary}\label{cor:commute-generic}
Let $b \in K$ be such that $b^{\tildefree} \in W$ and
$b^{\thefree}$ is algebraically independent over $K_{0}$.
Then, $\condition$ commutes.
\end{corollary}
\begin{proof}
Remember that we are assuming that $W$ is $K_{0}$-irreducible.
By Remark~\ref{rem:g-iteration}, for every
$\mu \in \Theta$ and $g, g' \in G_{\mu}$, $g(b^{\tildefree}) = g'(b^{\tildefree})$.
Since $W$ is irreducible and $b^{\tildefree}$ is generic in $W$ (over~$K_{0}$),
we have that $g = g'$.
\end{proof}

\subsection{The axioms}
\begin{definition}
The axioms of $\Tdbg$ are the axioms of $\Tdb$ plus the following axiom
scheme:
\begin{sentence}[(\texttt{$k$-Deep})]
Let $\condition = (W; p_{\pi}: \pi \in \thepivot)$ be a commutative \posture.%
\footnote{We are no longer assuming that $W$ is irreducible.}\linebreak[3]
Let $X \subseteq \K^{\tildefree}$ be $L(\K)$-definable, such that
$\Pi_{\thefree}(X)$ is large.
Then, there exists $a \in \K$ such that $a^{\tildefree} \in X$.
\end{sentence}
\end{definition}
Notice that the above is the analogue of the axiom scheme \ref{deep}. We don't
have an analogue for the axiom scheme \ref{wide}.
Notice also that the axiom scheme (\texttt{$k$-Deep}) is first-order expressible thanks to
Theorem~\ref{thm:local-global-commute}. 

\begin{thm}
\begin{enumerate}
\item $\Tdbg$ is a consistent and complete extension of $\Tdb$.  
\item If $T$ is model-complete, 
then $\Tdbg$ is an axiomatization for the model completion of $\Tdb$.
\item If $T$ eliminates quantifiers, then $\Tdbg$ eliminates quantifiers.
\item For every $\Ldeltabar$-formula $\alpha(\x)$ there exists an $L$-formula $\beta(\x)$
such that
\[
\Tdbg \models \forall \x \Pa{\alpha(\x) \leftrightarrow \beta(\x^{\Theta})}
\]
\item For every $\tuple{\K, \deltabar} \models \Tdbg$, for every $\av$ tuple in $\K$
and $B \subseteq \K$, the $\Ldeltabar$-type of $\av$ over $B$ is determined by the
$L$-type of $\av^{\Theta}$ over $B^{\Theta}$.
\end{enumerate}
\end{thm}

We assume that $T$ eliminates quantifiers.
We use the criterion in Proposition~\ref{prop:MC}(3) to show that $\Tdbg$ is the
model completion of $\Tdb$ and it eliminates quantifiers.
We will do it in two lemmas.

\begin{lemma}
Let $\tuple{A, \deltabar} \models \Tdb$.  
Let $\condition = (W; p_{\pi}: \pi \in \thepivot)$ be a commutative \posture
with parameters in~$A$, and
 $X \subseteq \K^{\tildefree}$ be an $L(A)$-definable set, such that
$\Pi_{\thefree}(X)$ is large.

Then, there exists $\tuple{B,  \deltabar} \supseteq \tuple{A, \deltabar}$ and $b \in
B$ such that $B \succeq A$,
$\tuple{B, \deltabar} \models \Tdb$, and $b^{\tildefree} \in X$.
\end{lemma}
\begin{proof}
Let $B \succeq A$ be $\card A^{+}$-saturated.
Let $\bv = (b_{\nu}: \nu \in \tildefree) \in X_{B}$ be such that 
$b_{\thefree} := \Pi_{\thefree}(\bv)$ 
is algebraically independent over~$A$.
Let $I \subset B$ be such that $I$ is disjoint from $b_{\thefree}$ and $D \coloneqq b_{\thefree} \cup I$ is a
transcendence basis of~$B$ over~$A$.
We extend each derivation $\delta \in \deltabar$ to $B$ in the following way.
It suffices to specify the value of $\delta$ on each $c \in D$.

If $c \in I$, we define $\delta c \coloneqq 0$.

If $c = b_{\mu}$ for some $\mu \in \thefree$, we define
$\delta b_{\mu} \coloneqq g_{\delta\mu}(\bv)$.

By definition, it is clear that $b_{0}^{\tildefree} = \bv \in X$.
Thus, it suffices to show that the extensions of $\deltabar$ commute on
all~$B$.
Again, it suffices to show that, for every $c \in D$ and every $\delta, \eps \in
\deltabar$,
\[
\delta \eps c= \eps \delta c.
\]
If $c \in I$, then both sides are equal to $0$, and we are done.

Suppose now that $c = b_{\mu}$ for some $\mu \in \thefree$, we have to show that
\begin{equation}
\label{eq:6}
\delta \eps b_{\mu}= \eps \delta b_{\mu}.
\end{equation}

Since the argument is delicate we prefer to give two different proofs with two different approaches.\\

\textsc{First proof)}
For this first proof, we replace $W$ with its $A$-irreducible component
containing $\bv$. Thus, we may assume that $W$ is $A$-irreducible.
Since $\bv$ is generic in $W$, the map $A(W) \to A(\bv)$, $h \mapsto h(\bv)$ is an
isomorphism of $A$-algebrae.
Via the above isomorphism, the derivation $R^{\delta}$ on $A(W)$ corresponds to the
derivation $\delta$ on $A(\bv)$.
The assumptions plus Theorem~\ref{thm:local-global-commute} and Remark~\ref{rem:R-comm} 
imply that $R^{\delta}$ and $R^{\eps}$ commute: therefore, also $\delta$ and $\eps$
commute on~$A(\bv)$.\\

\textsc{Second proof)}
Since $\bv \in W$, we have that $\delta b_{\pi} = f_{\delta, \pi}(\bv)$ for every $\pi \in
\thepivot$.
Thus, by definition, for every $\nu \in \tildefree$,
\[
\delta b_{\nu} = f_{\delta, \nu}(\bv).
\]

Denote $f \coloneqq f_{\eps \mu}$.
By definition, the LHS of \eqref{eq:6} is equal to
\begin{multline}
\delta g_{\eps \mu}(\bv) 
= \delta f(\bv) 
= f^{\delta}(\bv) + \sum_{\nu \in \tildefree} \frac{\partial f}{\partial \nu}(\bv) \cdot \delta b_{\nu} =\\
=  f^{\delta}(\bv) + \sum_{\nu \in \tildefree} \frac{\partial f}{\partial \nu}(\bv) \cdot f_{\delta, \nu}(\bv)
= f_{\delta, \eps \mu}(\bv) = g_{\delta, \eps \mu}(\bv).
\end{multline}
Similarly, the RHS of \eqref{eq:6} is equal to
$g_{\eps, \delta \mu}(\bv)$.
Finally, since $\condition$ commutes,
$g_{\delta, \eps \mu} = g_{\eps, \delta \mu}$, and we are done.
\end{proof}

\begin{lemma}\label{lem:Tdb-amalgamation}
Let $\tuple{B, \deltabar} \models \Tdb$,  $\tuple{C,\deltabar} \models \Tdbg$,
and $\tuple{A, \deltabar}$ be a common substructure,
such that $B$ and $C$ have the same $L(A)$-theory.
Let $\gamma(x)$ be a quantifier-free $\Ldb$-formula with parameters in $A$.
Let $b \in B$ be such that $\tuple{B, \deltabar} \models \gamma(b)$.
Then, there exists $c \in C$ such that $\tuple{C, \deltabar} \models \gamma(c)$.
\end{lemma}
\begin{proof}
\Wlog, we may assume that the only constants in the language $L$ are the
elements of $F$, 
the only function symbols are $+$ and $\cdot$, thus, an $L(\deltabar)$-substructure is a
differential subring containing~$F$.
Let $A'$ (resp., $A''$)
be the relative algebraic closure (as fields, or equivalently as $L$-structures)  of $A$ inside $B$ (resp., $C$).
Since $A$ has the same $L$-type in $B$ and $C$, there exists an isomorphism $\phi$
of $L$-structures between of $A'$ and $A''$ extending the identity on~$A$.
Moreover, any derivations on $A$ extends uniquely to $A'$ and $A''$: thus, 
$\phi$ is also an isomorphism of $\Ldeltabar$-structures.
Thus, \wloG we may assume that $A$ is relatively algebraically closed in $B$ and
in~$C$.

Let
\[
\thebounded \coloneqq \set {\mu \in \Theta: \mu b \in \acl(A b^{< \mu})}
\]
and $\thefree \coloneqq \Theta \setminus \thebounded$.
If $\thefree$ is empty, then $b \in \acl(A) = A$, and we are done.

Otherwise, let $\thepivot$ be the set of $\preceq$-minimal elements of~$\thebounded$.
For every $\pi \in \thepivot$, let $p_{\pi} \in A[x_{\pi}, x_{\thefree< \pi}]$ be such that
$p_{\pi}(\pi b, b^{\thefree < \pi}) =^{\pi} 0$.
Let $W$ be the $A$-irreducible component of 
\[
W_{0} \coloneqq \set{d_{\tildefree} \in B^{\tildefree}: \bigwedge_{\pi \in \thepivot}
p_{\pi}(d_{\tildefree}) =^{\pi} 0
}
\]
containing $b^{\tildefree}$.

\begin{claim}
$\condition \coloneqq (W; p_{\pi}: \pi \in \thepivot)$ is a commutative \posture
(over~$A$).
\end{claim}

Again, we prefer to give two different proofs.\\

\textsc{First proof)} Since $b^{\tildefree}$ is generic (over $A$)
in the $A$-irreducible variety~$W$, the map  $\Phi: A(W) \mapsto A(b^{\tildefree})$, $h \mapsto h(b^{\tildefree})$ is an
isomorphism of $A$-algebrae.
Via the above isomorphism, $R^{\delta}$ corresponds to the derivation $\delta$ on
$A(b^{\tildefree})= A(b^{\Theta})$: thus, the maps $R^{\delta}$ commute with each other.
Therefore, for every $\pi, \pi' \in \thepivot$ and $w,w' \in \Gamma$ with $[w] \pi = [w'] \pi'$,
\[
g_{w, \pi} = R^{w} g_{\pi} = \Phi((b^{\pi})^{[w]}) = 
\Phi((b^{\pi'})^{[w']}) = g_{w', \pi'}.
\]
\\

\textsc{Second proof)}
By Corollary~\ref{cor:commute-generic}.

\smallskip

Let $\beta(z_{\Theta})$ be a quantifier-free $L(A)$-formula such that
\[
T^{\deltabar} \cup \Diag_{L}(A) \vdash \gamma(x) \leftrightarrow \beta(x^{\Theta}).
\]

Let 
\[
X \coloneqq \set{x_{\tildefree} \in W: 
\dv \in W : B \models \beta \Pa{(g_{\mu}(\dv ))_{\mu \in \Theta}}
}.
\]
Since $b^{\tildefree} \in X$, and $X$ is $L(A)$-definable,
we have that $\Pi_{\thefree}(X)$ is large.
Thus, there exists $c \in C$ such that $c^{\tildefree} \in X$.
Since $c^{\tildefree} \in W$, by Remark~\ref{rem:g-iteration}
for every $\mu \in \Theta$ we have 
$c^{\mu} = g_{\mu}(c^{\tildefree})$.
Therefore, 
$\tuple{C, \deltabar} \models \beta( c^{\Theta})$.
\end{proof}

\subsubsection{Addendum}
Call a \posture $\condition =(W; p_{\pi}: \pi \in \thepivot)$ ``$K_0$-irreducible'' if $W$ is $K_{0}$-irreducible.
In the definition of a \posture we did not impose that it is
irreducible. However, by the proof of the Amalgamation
Lemma~\ref{lem:Tdb-amalgamation}, it seems that it would suffice to impose
in Axiom~(\texttt{$k$-Deep}) that only irreducible \postures
need to be satisfied.
The reason why we did not restrict ourselves to irreducible \postures is that we
don't know if we can impose irreducibility in a first-order way.
\begin{question}\label{Q:irreducible}
Let $(W_{i}: i \in I)$ be an $L$-definable family of varieties in $K^{n}$.
Is the set
\[
\set{i \in I: W_{i} \text{ is $K$-irreducible}} 
\]
definable?
\end{question}

\section{Stability and NIP}\label{sec:stable}

In this section we see that some of the model theoretic properties of T are inherited by $\Tdsg.$ In a following paper we will consider other properties. 
 We assume basic knowledge about stable and NIP theories: see \cite{Simon, Pillay}.

\begin{thm}
\begin{enumerate}
\item If $T$ is stable, then $\Tdsg$ is stable.
\item If $T$ is NIP, then $\Tdsg$ is NIP.
\end{enumerate}
\end{thm}
 
The above theorem follows immediately from the following one.
\begin{thm}
Let $U$ be an $L$-theory.
Let $\deltabar$ be a set of new \textbf{unary} function symbols.
Let $U'$ be an $\Ldb$-theory expanding~$U$.
Assume that, for every $\Ldb$-formula $\alpha(\x)$ there exists and $L$-formula
$\beta(\y)$ such that
\[
U' \models \forall \x\ \alpha(\x) \leftrightarrow \beta(\x^{\Gamma}),
\]
where $\x^{\Gamma}$ is the set of $\deltabar$-terms in the variables $\x$.

Then, for every $(M, \deltabar) \models U'$ and every $\av$ tuple in $M$ and $B$
subset of $M$, the $\Ldb$-type of $\av$ over $B$ is uniquely determined by the $L$-typle
of $\av^{\Gamma}$ over $B^{\Gamma}$.

Moreover, 
\begin{enumerate}
\item If $U$ is stable, then $U'$ is stable.
\item If $U$ NIP, then $U'$ is NIP.
\end{enumerate}
\end{thm}

\begin{proof}
The results follow easily by applying the following criteria.

1) \cite[Thm II. 2.13]{Shelah} A theory $U$ is stable iff, for every subset $A$ of a model $M$ of $U$, and for
every sequence $(\av_{n})_{n \in \N}$ of tuples in $M$, if $(\av_{n})_{n \in \N}$ is
an indiscernible sequence, then it is totally indiscernible.

2) 
\cite[Proposition 2.8]{Simon} A theory $U$ is NIP iff, 
 for every formula $\phi(\x; \y)$ and for any
indiscernible sequence $(\av_{i} : i \in I)$ and tuple $\bv$, 
there is some end segment
$I_{0} \subseteq I$ such that $\phi(a_{i} ; b)$ is ``constant'' on $I_{0}$: that is, 
either for every $i \in I_{0}$   $\phi(\av_{i} ; \bv)$ holds, 
or for every $i \in I_{0}$   $\neg \phi(\av_{i} ; \bv)$ holds.
\end{proof}

\section{Pierce-Pillay axioms}\label{sec:PP}
We give now an extra axiomatization for $\Tdg$, in the ``geometric'' style of Pierce and Pillay~\cite{PP}.
We won't use this axiomatization, but it may be of interest.

Let $\tuple{\K, \delta} \models \Td$.


Let $W \subseteq \K^{n}$ be an algebraic variety defined over~$\K$.
We define the twisted tangent bundle $\tau W$ of $W$  \wrt $\delta$ in the same way as in 
\cite{PP} (see also \cite{Moosa:22}, where it is called ``prolongation'').

Let $\x \coloneqq \tuple{x_{1}, \dotsc, x_{n}}$.
Let $\K^{*} \succeq \K$ and $\av \coloneqq \tuple{a_{1}, \dotsc, a_{n}}\in
(\K^{*})^{n}$.
We define
\[\begin{aligned}
I(\av/\K) \coloneqq& \set{p \in \K[\x]: p(\av) = 0}\\
I(W/\K) \coloneqq& \set{p \in \K[\x]: \forall \cv \in W\ p(\cv) = 0 }.
\end{aligned}\]

Let $\bar p = (p_{1}, \dotsc, p_{\ell}) \in \K[\x]^{\ell}$.
Define 
\[
V_{\K}(\bar p) \coloneqq \set{\cv \in \K^{n}: p_{i}(\cv) = 0, i = 1, \dotsc, \ell},
\]
and $(\bar p)_{\K}$ to be the ideal of $\K[\x]$ generated by $p_{1}, \dotsc, p_{\ell}$.

\begin{definition}\label{def:torsor}
Assume that 
$I(W/\K) = (p_{1}, \dotsc, p_{\ell})_{\K}$.
The twisted tangent bundle of~$W$ (\wrt~$\delta$) is the algebraic variety $\tau^{\delta} W \subseteq \K^{n} \times \K^{n}$
\[
\tau^{\delta} W  \coloneqq \set{\tuple{\x, \y} \in \K^{n} \times \K^{n}:
p_{i}^{[\delta]}(\x,\y) = 0, i=1, \dotsc, n}
\subseteq \K^{n} \times \K^{n}, 
\]
where $p^{[\delta]}$ was introduced in Definition~\ref{def:delta-poly}.
\end{definition}
Notice that the definition of $\tau^{\delta} W$ does not depend on the choice of polynomials
$\bar p$ such that $I(W/\K) = (\bar p)_{\K}$.
Notice also that, when $\delta = 0$, the twisted tangent bundle $\tau^{0}W$ coincides with the tangent bundle.

The importance in this context of the twisted tangent bundle is due to the following two facts:
\begin{remark}
If $\av \in W$, then $\tuple{\av, \delta \av} \in \tau^{\delta} W$.
\end{remark}
\begin{fact}\label{fact:torsor}
Let $L \supset \K$ be a field.
Let $\bv \in W$ be (as interpreted in~$L$) such that $\bv$ is generic in $W$ over~$\K$ (that is,
$\rk(\bv / \K) = \dim(W)$).
Let $\cv \in L^{n}$ be such that $\tuple{\bv,\cv} \in \tau^{\delta} W$.

Then, there exists a derivation $\eps$ on $L$ extending $\delta$
and such that $\eps \bv = \cv$.
\end{fact}
\begin{proof}
It is a known result: see  \cite[Ch.~II, \S17, Thm.~39]{ZS1}.
See also \cite[Thm.~VIII.5.1]{Lang}, \cite{PP},
and \cite[Lemma~1.1]{GR}.
\end{proof}

We want to write an axiom scheme generalizing Pierce-Pillay to $\Tdg$.

An idea would be to use the following:
\begin{customenum}
\setcounter{customenumi}{3}
\item
\label{ax:PPw} Let  $W \subseteq \K^{n}$ be an algebraic variety which is defined
over $\K$ and $\K$-irreducible. Let $U \subseteq \tau^{\delta} W$ be an $L(\K)$-definable set, such
that the projection of $U$ over $W$ is large in~$W$ (\ie, of the same dimension
as~$W$).
Then, there exists $\av \in W$ such that $\tuple{\av, \delta \av} \in U$.
\end{customenum}

However, there is an issue with the above axiom scheme: we don't know how to
express it in a first order way!
The reason is the following: give a definable family of tuples of polynomials
$(\bar p_{i}: i \in I)$, while each $\tau^{\delta}(V_{\K}(\bar p_{i}))$ is definable, we
do not know whether the family $\Pa{\tau^{\delta}(V_{\K}(\bar p_{i})): i \in I}$ is
definable.
We leave it as an open problem, and we will use a different axiom scheme.

\begin{question}
Let $\Pa{\bar p_{i}: i \in I}$ be a definable family of tuples of polynomials.
Is there a definable family $\Pa{\bar q_{i}: i \in I}$ of tuples of polynomials,
such that $I(V_{\K}(\bar p_{i})/\K) = (\bar q_{i})_{\K}$ for every $i \in I$?
\end{question}
The above question is related to Question~\ref{Q:irreducible}:
notice that 
``$W$~is $\K$-irreducible'' is equivalent to ``$I(W/\K)$ is prime'', and
the latter, by \cite{DS} (see also \cite{Schoutens}),   is a definable property
of the parameters of the formula defining~$I(W/\K)$.

\medskip

We need some additional definitions and results 
before introducing the true axiom scheme.
Fix $\pv \in \K[\x]^{\ell}$, and let $W \coloneqq V_{\K}(\pv)$.
Given $\av \in W$, the twisted tangent space of $\pv$ at $\av$ is
\[
\tau^{\delta}_{\av}(\pv) \coloneqq \set{\y \in \K^{n}: 
p_{i}^{[\delta]}(\av, \y) = 0: i = 1, \dotsc, \ell}.
\]
Moreover, $\tau^{0}(\pv)$ is the usual tangent space at $\av$ of $V_{\K}(\pv)$.

\begin{remark}\label{rem:tangent-dimension}
Let $\av \in ({K^{*}})^{n}$, we define $J \coloneqq I(\av/\K)$ and $F \coloneqq \K(\av) \subseteq \K^{*}$.
Let $J'$ be the ideal of $F[\x]$ generated by~$J$ and
let $S \coloneqq \tau^{0}_{\av}J'$ be the tangent space at~$\av$ 
(as an $F$-vector space).
Then, the dimension of $S$ as an $F$-vector space is equal to $\rk(\av/\K)$.
Indeed, let $S'$ be the set of derivations on $F$ which are $0$ on~$\K$: then,
by Fact~\ref{fact:torsor}, $S$ and $S'$ are isomorphic as $F$-vector spaces.
By \cite[Ch.~II, \S17, Thm.~41]{ZS1}, the dimension of $S'$ as $F$-vector space
is equal to $\rk(\av/\K)$.
\end{remark}

Let $W \coloneqq V_{\K}(\pv)$;
the following lemma shows that, under suitable conditions, we can replace
$I(W/\K)$ with $I(\pv/\K)$. We need to introduce some notations:

Let $d_{0} \coloneqq \dim(W)$ and for every $d \in \N$, we define
\[
\Reg^{d}(\pv) \coloneqq \set{\cv \in W: \dim(\tau^{0}_{\av}(\pv)) = d}
\]
and $\Reg(\pv) \coloneqq \Reg^{d_{0}}(\pv)$.

\begin{lemma}\label{lem:PP-generic}
Let $\av \in ({\K^{*}})^{n}$ be such that $\av \in \Reg(\pv)$ and $rk(\av/\K) \geq d_{0}$.
Let $J \coloneqq I(\av/\K)$ and $W \coloneqq V_{\K}(\pv)$.

Then,
\begin{enumerate}
\item $\dim(\Reg(\pv)) = \rk(\av/\K) = d_{0}$;
\item $\tau^{0}_{\av}(\pv) = \tau^{0}_{\av}(J)$;
\item $\tau^{\delta}_{\av}(\pv) = \tau^{\delta}_{\av}(J)$;
\item for every $\bv \in ({\K^{*}})^{n}$ such that
$\tuple{\av, \bv} \in \tau^{\delta}(\pv)$ there exists $\delta'$ derivation on $\K^{*}$ extending
$\delta$ and such that $\delta'(\av) = \bv$.
\end{enumerate}
\end{lemma}
\begin{proof}
1) It is clear:
\[
d_{0} \leq \rk(\av/\K) \leq \dim(\Reg(\pv)) \leq \dim(W) = d_{0}.
\]

2) Since $\pv(\av) = 0$, we have $(\pv)_{\K} \subseteq J$, and therefore
\[
V_{\K}(J) \subseteq V_{\K}(\pv) = W.
\]
Thus, $\dim(V_{\K}(J)) \leq \dim W = d_{0}$.
Therefore, $\av$ is also a generic point of $V_{\K}(J)$.
So, by Remark~\ref{rem:tangent-dimension},
\[
dim(\tau^{0}_{\av}(J)) = \rk(\av/\K) = d_{0}.
\]
Therefore, the vector space
 $\tau^{0}_{\av}(\pv)$ contains $\tau^{0}_{\av}(J)$ and has the same dimension
 $d_{0}$:
thus, they are equal.

3) Notice that
 $\tau^{0}_{\av}(\pv)$ and $\tau^{\delta}_{\av}(\pv)$ are vector spaces of the same
 dimension, and the same happens for $\tau^{\delta}_{\av}(J)$.
Moreover, $\tau^{\delta}_{\av}(\pv)$ contains $\tau^{\delta}_{\av}(J)$ and has the same
dimension, and therefore $\tau^{\delta}_{\av}(\pv) = \tau^{\delta}_{\av}(J)$.

4) It follows from Fact~\ref{fact:torsor}.
\end{proof}

Given $m, n, d \in \N$, let
\[
(\pv_{m,n,d}(\x, \av): \av \in \K^{\ell})
\]
be a parametrization (definable in the language of rings) of all $m$-tuples of
polynomials  in $\K[x_{1}, \dotsc, x_{n}]$ of degree at most~$d$.
We will write $\pv$ instead of $\pv_{m,n,d}$.

For each $\av \in \K^{\ell}$, let $W_{\av} \coloneqq V_{\K}(\pv(\x, \av))$ and
$U_{\av} \coloneqq \Reg(\pv(\x,\av) )\subseteq W_{\av}$.
Let $\Pi_{n}: \K^{n} \times \K^{n} \to \K^{n}$ be the canonical projection onto the first
$n$ coordinates.

We can finally write the axiom scheme.

\begin{customenum}
\setcounter{customenumi}{2}
\item
\label{ax:PP} Let $m,n,d \in \N$ and $\pv \coloneqq \pv_{m,n,d}(\x,\y)$.
Let $\av \in \K^{\ell}$ and $X \subseteq \tau^{\delta}(\pv(\x, \av))$ be $L(\K)$-definable.
Assume  that
$\Pi_{n}(X) \subseteq U_{\av}$ and $\dim(\Pi_{n}(X)) = \dim(W_{\av})$.
Then, there exists $\cv \in \K^{n}$ such that $\tuple{\cv, \delta \cv} \in X$.
\end{customenum}

\begin{thm}
$\TdgPP \coloneqq \Td \cup \text{\ref{ax:PP}}$ is an axiomatization of $\Tdg$.   
\end{thm}
\begin{proof}
Since we can take $\pv$ to be the empty tuple,
and therefore $W_{\empty} = \K^{n}$, it is clear that \ref{ax:PP} implies
\ref{wide}.

We have to prove the opposite.
Since $\Tdg$ is complete, it suffices to show that $\TdgPP$ is consistent. 
\Wlog, we may assume that $T$ has elimination of quantifiers.
To show that $\TdgPP$ is consistent, it suffices to prove the following
\begin{claim}
Let $m,n,d, \pv, \av, X$ be as in \ref{ax:PP}.
Then there exists $K^{*} \succeq \K$ and $\bv \in ({\K^{*}})^{n}$ such that
$\tuple{\bv, \delta \bv} \in X$.
\end{claim}
Let $K^{*} \succ \K$ be sufficiently saturated and $\bv \in {\K^{*}}^{n}$ such that
$\bv \in \Pi_{n}(X)$ and $\rk(\bv/K) = \dim(\Pi_{n}(X)) = \dim(W_{\av})$.
Let $\cv \in {\K^{*}}^{n}$ be such that $\tuple{\bv, \cv} \in X$.
By Lemma~\ref{lem:PP-generic} there exists $\delta'$ derivation on $\K^{*}$ extending
$\delta$ and such that $\delta(\bv) = \cv$.
\end{proof}

Giving the analogue axiomatization for $\Tdng$ is not difficult and the reader can provide the details.

On the other hand, we won't try to give a similar axiomatization for~$\Tdbg$, 
since already when $T= ACF$ it is an arduous task: see 
\cites{Pierce14,  Pierce:2, LS}.

\section{Conjectures and open problems}

We conclude the paper with a list of open problems, remarks and some ideas.

\subsection{Elimination of imaginaries}

\begin{conjecture}
$\Tdsg$ has  elimination of imaginaries modulo $T^{eq}$.
\end{conjecture}
A few particular cases are known, when $\Tdsg$ is one of the following:
\begin{itemize}
\item $\DCFOm$: see \cite{mcgrail};
\item $\rm{RCF}$ or certain theories of Henselian valued fields, endowed with $m$ commuting generic derivations: see \cite{FK, KP} for a
proof based on M.~Tressl's idea; 
see also \cites{bkp, Point, Silvain} for different proofs;
\item $\DCFOmnc$ (see \cite{MS}).
\end{itemize}
We have also established the validity of the above conjecture for certain topological structures. 
Drawing upon established techniques, it is probable that the conjecture can be
proven for $T$ simple (as proved in \cites{MS, Mohamed}). However, for the
general case, we believe that novel approaches are required (although some
progress has been made in \cite{bkp}).


\subsection{Definable types}
Let $\tuple{\K, \deltabar}\models \Tdsg$.
Given a type $p \in S_{\Ldelta}^{n}(\K)$, let $\av$ be a realization of $p$; we
define $\tilde p \in S_{L}^{n\times \Gamma}(\K)$ as the $L$-type of $\av^{\Gamma}$ over~$\K$.
\begin{open problem}
Is it true that $p$ is definable iff $\tilde p$ is definable?
We conjecture that it is true when $\Tdsg = \Tdbg$.
\end{open problem}

\subsection{Zariski closure}
Given $X \subseteq \K^{n}$, denote by $X^{Zar}$ be the Zariski closure of $X$.
\begin{open problem}[See \cite{FLL}]
1) Let $\Pa{X_{i}: i \in I}$ be an $L$-definable family of subsets of $\K^{n}$.
Is $\Pa{X_{i}^{Zar}: i \in I}$ also $L$-definable?

2) Assume that 1) holds for $\K$.
Let $\tuple{\K, \deltabar} \models \Tdsg$.
Let $\Pa{X_{i}: i \in I}$ be an $\Ldelta$-definable family of subsets of $\K^{n}$.
Is $\Pa{X_{i}^{Zar}: i \in I}$ also $\Ldelta$-definable?
\end{open problem}

\subsection{Monoid actions}
Let $\Lambda$ be a monoid generated by a $k$-tuple $\deltabar$: we consider $\Lambda$ as a
quotient of the free monoid~$\Gamma$.
We can consider actions of $\Lambda$ on models of $T$ such that each $\delta_{i}$ is a
derivation: we have a corresponding theory $T^{\Lambda}$ whose language is $\Ldelta$
and with axioms given by $T$, the conditions that each $\delta_{i}$ is a derivation,
and, for every $\gamma, \gamma' \in \Gamma$ which induce the same element of $\Lambda$, the axiom
 $\forall x\, \gamma x = \gamma'x$.
\begin{open problem}
Under which conditions on $\Lambda$ the theory $T^{\Lambda}$ has a model completion?
\end{open problem}

\begin{conjecture}
Let $\Gamma_{\ell}$ be the free monoid in $\ell$ generators, 
and $\Theta_{k}$ be the free commutative monoid in $k$ generators.
Then, for $\Lambda$  equal either to
$\Gamma_{\ell} \times \Theta_{k}$ or to $\Gamma_{\ell} * \theta_{k}$, $T^{\Lambda}$ has a model completion 
(where $\times$ is the Cartesian product, and $*$ is the free product).
More generally, for $\Gamma$ equal to a combination of free and Cartesian products of
finitely  many copies of $\N$, $T^{\Lambda}$ has a model completion.
\end{conjecture}
As a consequence of \cite{LeonM:16}, when $T = ACF_{0}$ and $\Lambda = \Gamma_{\ell}\times \Theta_{k}$,
 $T^{\Lambda}$ has a model companion.

Maybe the following conditions on $\Lambda$ suffice for $T^{\Lambda}$ to have a model completion:\\
Let $\preceq$ be the canonical  quasi ordering on $\Lambda$ given by $\alpha \preceq \beta \alpha$ for every 
$\alpha,\beta \in \Lambda$; we assume that:
\begin{itemize}
\item  $\preceq$ is a well-founded partial ordering;
\item for every $\lambda \in \Lambda$, the set $\set{\alpha \in \Lambda: \alpha \preceq \lambda}$ is finite;
\item for every $\alpha, \beta \in \Lambda$, if they have an upper bound, then they have a least
upper bound;
\item let $X \subset \Lambda$ be finite; assume that $X$ is $\preceq$-initial in $\Lambda$; then,
$\Lambda \setminus X$ has finitely many $\preceq$-minimal elements;
\item if $\alpha_{1} \delta_{1} = \alpha_{2} \delta_{2}$ for some $\alpha_{i} \in \Lambda$ and $\delta_{i} \in
\deltabar$,
then $\delta_{1}$ and $\delta_{2}$ commute with each other; moreover, there exists
$\beta \in \Lambda$ such that $\alpha_{1} = \delta_{2} \beta$ and $\alpha_{2} = \delta_{1} \beta$.
\end{itemize}

\section*{Aknwoledgements}
The authors thank Noa Lavi, Giorgio Ottaviani, Françoise Point, 
Silvain Rideau-Kikuchi,
Omar Le\'on Sanchez, and Marcus Tressl for the interesting discussions on the topic.

\section*{Funding}
Both authors are members of the ``National Group for Algebraic and Geometric Structures,
 and their Applications'' (GNSAGA - INdAM). 
This research is part of the project  PRIN 2022 ``Models, sets and classification''.

\printbibliography

\end{document}